\pdfoutput=1
\documentclass[journal]{IEEEtran}

\IEEEoverridecommandlockouts   

\usepackage{graphics} 
\usepackage{epsfig} 
\usepackage{times} 
\usepackage{algorithmic}
\usepackage{array}
\usepackage{amsmath} 
\usepackage{amsthm}
\usepackage{amssymb}  
 \usepackage{graphicx}
\usepackage{epstopdf}
\usepackage{cases}
\usepackage{subfig}
 \usepackage{dblfloatfix}
\usepackage{url}
\usepackage{color}
\usepackage{bm}
\usepackage{booktabs} 
\usepackage{stmaryrd}
\usepackage{cite}
\usepackage{ifpdf}
\newtheorem{Theorem 1}{Theorem}
\newtheorem{Theorem 2}[Theorem 1]{Theorem}
\newtheorem{Definition 1}{Definition}
\newtheorem{Definition 2}[Definition 1]{Definition}
\newtheorem{Remark 1}{Remark}
\newtheorem{Remark 2}[Remark 1]{Remark}
\newtheorem{Lemma 1}{Lemma}
\newtheorem{Lemma 2}[Lemma 1]{Lemma}
\newtheorem{Assumption 1}{Assumption}
\newtheorem{Corollary 1}{Corollary}
\newtheorem*{Update rule}{Update rule for $b_{i \shortrightarrow  j}^t$}
\newtheorem*{Condition A}{Condition A}
\newtheorem*{Condition B}{Condition B}
\usepackage{multirow}

\title{\LARGE \bf
	ADMM Based Privacy-preserving Decentralized Optimization }

\author{Chunlei Zhang, Muaz Ahmad, Yongqiang Wang
	\thanks{*The work was supported in part by the National Science Foundation under Grant 1738902.}
	\thanks{Chunlei Zhang, Muaz Ahmad and Yongqiang Wang are with the
		department of Electrical and Computer Engineering, Clemson
		University, Clemson, SC 29634, USA
		{\tt\small \{chunlez,muaza,yongqiw\}@clemson.edu}}%
}

\begin{document}

\maketitle
\thispagestyle{empty}
\pagestyle{empty}


\begin{abstract}
    Privacy preservation is addressed for decentralized optimization, where $N$ agents cooperatively minimize  the sum of $N$  convex functions private to these individual agents.
    In most existing decentralized optimization approaches,  participating agents exchange and disclose states explicitly, which may not be desirable when the states contain sensitive information of individual agents. The problem is more acute when adversaries exist which try  to steal information from other participating agents. To address this issue, we propose a privacy-preserving decentralized optimization approach based on ADMM and partially homomorphic cryptography. To our knowledge, this is the first time that cryptographic techniques are incorporated in a fully decentralized setting to enable privacy preservation in decentralized optimization in the absence of any third party or aggregator. To facilitate the incorporation of encryption in a fully decentralized manner, we introduce a new ADMM which allows time-varying penalty matrices and rigorously prove that it has a convergence rate of $O(1/t)$. Numerical and experimental results confirm the effectiveness and low computational complexity of the  proposed approach. 
\end{abstract}

\section{Introduction}
In recent years, decentralized optimization has been playing key roles in applications as diverse as rendezvous in multi-agent systems \cite{lin2004multi}, spectrum sensing in cognitive networks \cite{zeng2011distributed}, support vector machine in machine learning \cite{cortes1995support}, online learning \cite{yan2013distributed}, classification \cite{zhang2017dynamic}, data regression in statistics \cite{mateos2010distributed}, source localization in sensor networks \cite{zhang2017distributed}, and monitoring of smart grids \cite{nabavi2015distributed}. 
In these applications, the problem can be formulated in the following general form, in which $N$ agents cooperatively solve an unconstrained optimization problem:
\begin{equation}\label{eq:global function}
\begin{aligned}
\mathop {\min }\limits_{\bm{y}} \qquad \sum\limits_{i=1}^{N}f_i(\bm{y}), 
\end{aligned}
\end{equation}
where variable $\bm{y}\in \mathbb{R}^D$ is common to all agents, function $f_i: \mathbb{R}^D\to\mathbb{R}$ is the local objective function of agent $i$. 
 
Typical decentralized solutions to the optimization problem \eqref{eq:global function} include distributed (sub)gradient based algorithms \cite{nedic2009distributed}, augmented Lagrangian methods (ALM) \cite{he2016proximal}, and the alternating direction method of multipliers (ADMM) as well as its variants \cite{boyd2011distributed,he2016proximal,ling2013decentralized,ling2016weighted}, etc. In (sub)gradient based solutions, (sub)gradient computations and averaging among neighbors are conducted iteratively to achieve convergence to the minimum. In augmented Lagrangian  and ADMM based solutions, iterative Lagrangian  minimization is employed, which, coupled with dual variable update, guarantees that all agents agree on the same minimization solution.

However, most of the aforementioned decentralized approaches  require agents to exchange and disclose their states explicitly to neighboring agents in every iteration \cite{nedic2009distributed,boyd2011distributed,ling2016weighted,ling2013decentralized,he2016proximal}. This brings about serious privacy concerns in many practical applications  \cite{lagendijk2013encrypted}.  For example, in projection based source localization, intermediate states are positions of points lying on the circles centered at individual nodes' positions \cite{shi2008distributed}, and thus a node may infer the exact position of a neighboring node using three intermediate states, which is undesirable when agents want to keep their position  private  \cite{al2017proloc}.
In the rendezvous problem where a group of individuals want to meet at an agreed time and place \cite{lin2004multi}, exchanging explicit states may leak their initial locations which may need to be kept secret instead \cite{mo2017privacy}. Other examples include the agreement problem \cite{degroot1974reaching}, where a group of individuals want to reach consensus on a subject without leaking their individual opinions to others  \cite{mo2017privacy}, and the regression problem \cite{mateos2010distributed}, where individual agent's training data may contain sensitive information (e.g., salary, medical record) and should be kept private.  In addition, exchanging explicit	states without encryption is susceptible to eavesdroppers which try to intercept and steal information from exchanged messages.

To enable privacy preservation in decentralized optimization, one commonly used approach is differential privacy \cite{huang2015differentially,han2017differentially,nozari2016differentially}, which adds carefully-designed noise to exchanged states or objective functions to cover sensitive information. However, the added noise also unavoidably compromises the accuracy of optimization results, leading to a trade-off between privacy and accuracy \cite{huang2015differentially,han2017differentially,nozari2016differentially}. In fact, as indicated in  \cite{nozari2016differentially}, even when no noise perturbation is added, differential-privacy based approaches may fail to converge to the accurate optimal solution. It is worth noting that  although some differential-privacy based optimization approaches can  converge  to the optimal solution in the mean-square sense  with the assistance of a third party such as a cloud (e.g., \cite{hale2015differentially}, \cite{hale2017cloud}), those results are not applicable to the completely decentralized setting discussed here where no third parties or aggregators exist. Observability-based design has been proposed for privacy preservation in linear multi-agent networks \cite{alaeddini2017adaptive,pequito2014design}. By properly designing the weights for the communication graph, agents' information will not be revealed to non-neighboring agents. However, this approach cannot protect the privacy of the direct neighbors of compromised agents and  it is susceptible to external eavesdroppers. Another approach to enabling privacy preservation is encryption. However, despite successful applications in cloud based control and optimization  \cite{xu2015secure,freris2016distributed,shoukry2016privacy,wang2011secure},  conventional cryptographic techniques cannot be applied directly in a completely {\it decentralized} setting without the assistance of aggregators/third parties (note that traditional secure  multi-party computation schemes like fully homomorphic encryption \cite{gentry2009fully} and Yao's garbled circuit \cite{yao1982protocols} are computationally too  heavy to be practical for real-time optimization \cite{lagendijk2013encrypted}).  Other privacy-preserving optimization approaches include \cite{mangasarian2012privacy,gade2017private} which protect privacy via  perturbing problems or states. Recently, by using linear dynamical systems theory to facilitate cryptographic design, we proposed a privacy-preserving decentralized linear consensus approach \cite{ruan2017secure}. However, to the best of our knowledge, results are still lacking on cryptography based approaches that can enable privacy for the more complicated general nonlinear optimization problem like \eqref{eq:global function} in completely {\it decentralized} setting without any aggregator/third party.

Given that the compromised accuracy of differential-privacy based approaches makes them inappropriate for applications where both accuracy and privacy are of primary concern (for example, when dealing with medical treatment data, even minimal noise may interfere with users' learning \cite{jo2012safe,machanavajjhala2011personalized}), we propose a new privacy-preserving decentralized optimization approach based on ADMM and partially homomorphic cryptography. We used ADMM because it has several  advantages. First, ADMM has a fast convergence speed in both primal and dual iterations \cite{ling2016weighted}. By incorporating a quadratic regularization term, ADMM has been shown to be able to obtain satisfactory convergence speed even in ill-conditioned dual functions \cite{ling2013decentralized}. Secondly, the convergence of ADMM has been established for non-convex and non-smooth objective function \cite{wang2015global}. Moreover, from the implementation point of view, not only is ADMM easy to parallelize and implement, but it is also robust to noise and computation errors \cite{simonetto2014distributed}. 

It is worth noting that privacy has different meaning under different settings. For example, in the distributed optimization literature, privacy has been defined as the non-disclosure of  agents' states\cite{hale2015differentially}, objective functions or subgradients \cite{nozari2016differentially,yan2013distributed,lou2017privacy}. In this paper, we define privacy as preserving the confidentiality of agents' intermediate states,  gradients of objective functions, and objective functions.  We protect the privacy of objective functions  through protecting intermediate states. In fact, if left unprotected,  intermediate states could be used by an adversary to infer the gradients or even objective functions of other nodes through, e.g., data mining techniques. For example, in the regression problem in  \cite{mateos2010distributed}, the objective functions take the form  $f_i(\bm{y})=\frac{1}{2}\parallel\bm{s}_i-B_i\bm{y}\parallel_2^2$, in which $s_i$ and $B_i$ are raw data containing sensitive information such as salary and medical record. When the subgradient method in  \cite{nedic2009distributed} is used to solve the optimization problem $\mathop {\min }\limits_{\bm{y}} \quad \sum\limits_{i=1}^{N}f_i(\bm{y})$, agent $i$ updates its intermediate states in the following way:
	$$ \bm{y}_i^{k+1}=\sum_{j=1}^{N} a_{ij} \bm{y}_j^{k}+\alpha_k \triangledown f_i(\bm{y}_i^k) $$
	where $a_{ij}$ are weights, $\alpha_k$ is the stepsize, and $\triangledown f_i(\bm{y})=B_i^TB_i\bm{y}-B_i^T\bm{s}_i$ is the gradient. In  this case, an adversary can infer  $\triangledown f_i(\bm{y}_i^k)$ based on exchanged intermediate states $\bm{y}_i$ if the weights $a_{ij}$ and stepsize  $\alpha_i$ are publicly known.  We consider two adversaries:	{\it Honest-but-curious adversaries} are agents who follow all protocol steps correctly but are curious and collect all intermediate and input/output data in an attempt to learn some information about other participating agents \cite{li2010secure}. {\it External eavesdroppers} are adversaries who steal information through wiretapping all communication channels and intercepting exchanged  messages between agents. Protecting agents' intermediate states can avoid eavesdroppers from inferring any information in optimization.
	
\emph{Contributions:} The main contribution of this paper is a privacy-preserving decentralized optimization approach based on ADMM and partially homomorphic cryptography. To our knowledge, this is the first time that cryptographic techniques are incorporated in a fully decentralized setting to enable privacy preservation in decentralized optimization without the assistance of any third party or aggregator. To facilitate the incorporation of homomorphic encryption in ADMM in a fully decentralized manner, we also propose a new ADMM which allows time-varying penalty matrices and rigorously characterize its convergence rate of $O(1/t)$. It is worth noting that the $O(1/t)$ convergence rate requires the subproblems (primary update) to be efficiently computed. When the subproblems are difficult, some suboptimal solution (obtained by, e.g., the approach in \cite{barzilai1988two}) can be used to approximately solve the subproblem.  In contrast to differential-privacy based optimization approaches \cite{huang2015differentially,han2017differentially,hale2015differentially,nozari2016differentially},  our approach can enable privacy preservation without sacrificing accuracy.  Moreover, our approach does not require strong convexity, Lipschitz or twice continuous differentiability in the objective function. (It is worth noting that  \cite{hale2015differentially} requires a trusted cloud and hence is not applicable to the completely decentralized setting discussed here.) Different from the privacy-preserving optimization approach in \cite{lou2017privacy} which only protects the privacy of gradients, our approach preserves the privacy of both intermediate states and gradients. In addition, \cite{lou2017privacy} assumes that an adversary does not have access to the adjacency matrix of the network graph while our approach does not need this assumption.

\section{A New ADMM with Time-varying Penalty Matrices}
In this section, we propose a new ADMM with time-varying penalty matrices for \eqref{eq:global function}, which is key for enabling the incorporation of partially homomorphic cryptography in a completely decentralized optimization problem for privacy protection. 
\subsection{Problem Formulation}
We assume that each $f_i$ in \eqref{eq:global function} is private and only known to agent $i$, and all $N$ agents form a bidirectional connected network. Using the graph theory \cite{bondy1976graph}, we represent the communication pattern
of a multi-agent network by a graph $G=\{V,E\}$, where  $V$ denotes the set of agents and $E$ denotes the set of
communication links (undirected edges) between agents. Denote the total number of
communication links in $E$ as $|E|$. If there exists a communication link between agents $i$ and $j$, we say that agent $j$ is a neighbor of $i$ (agent $i$ is a neighbor of $j$ as well) and denote the communication link as
$e_{i,j}\in E$ if $i<j$ is true or $e_{j,i}\in E$ if $ i>j$ is true.  Moreover, we denote the set of all neighboring agents of $i$ as $\mathcal{N}_i$ (we consider agent $i$ to be a neighbor of itself in this paper, i.e., $i\in\mathcal{N}_i$, but $e_{i,i}\notin E$). 

\subsection{Proximal Jacobian ADMM}

To solve \eqref{eq:global function} in a decentralized manner, we reformulate \eqref{eq:global function} as follows (which avoids using dummy variables in conventional ADMM \cite{shi2014linear}):
  \begin{equation}\label{eq:conventional admm form}
  \begin{aligned}
  &\mathop {\min }\limits_{\bm{x}_{i}\in \mathbb{R}^D,\,i\in\{1,2,...,N\}} \qquad \sum\limits_{i=1}^{N}f_i(\bm{x}_{i})  \\
  &\textrm{subject to} \qquad \bm{x}_{i}=\bm{x}_{j}, \quad \forall e_{i,j}\in E,  
  \end{aligned}
  \end{equation}
  where $\bm{x}_i$ is considered as a copy of $\bm{x}$ and belongs to agent $i$. To solve \eqref{eq:conventional admm form}, each agent first exchanges its current state $\bm{x}_i$ 
 with its neighbors. Then it carries out local computations based on its private local objective function $f_i$ and the received state information from neighbors to update its state.  Iterating these computations will make every agent reach consensus on a solution that is optimal to \eqref{eq:global function} when \eqref{eq:global function} is convex. Detailed implementation of the ADMM algorithm based on Jacobian update is elaborated as follows \cite{deng2017parallel}:
 \begin{numcases}{}
\bm{x}_i^{t+1}=\underset{\bm{x}_i}{\operatorname{argmin}}\mathcal{L}(\bm{x}_1^t,\bm{x}_2^t,...,\bm{x}_{i-1}^t,\bm{x}_i,\bm{x}_{i+1}^t,...,\bm{x}_N^t,\bm{\lambda}^t) \nonumber \label{eq:conventional admm 1}\\
  \qquad\qquad\qquad\quad+\frac{\gamma_i}{2}\parallel \bm{x}_i-\bm{x}_i^t \parallel^2, \\
  \bm{\lambda}_{i,j}^{t+1}=\lambda_{i,j}^{t}+\rho(\bm{x}_i^{t+1}-\bm{x}_j^{t+1}), \quad\forall j\in\mathcal{N}_i \label{eq:conventional admm 2}
 \end{numcases}
 for $i=1,2,...,N$. Here, $t$ is the iteration index, $\gamma_i>0$ $(i=1,2,...,N)$ are proximal coefficients, and $\mathcal{L}$ is the augmented Lagrangian function
 \begin{equation} \label{eq:conventional augemented lagrangian function}
 \begin{aligned}
 \lefteqn{\mathcal{L}(\bm{x},\bm{\lambda})=\sum\limits_{i=1}^{N}f_i(\bm{x}_{i}) }\\
 & \qquad+ \sum\limits_{e_{i,j}\in E}(\bm{\lambda}_{i,j}^{T}(\bm{x}_i-\bm{x}_{j})+\frac{\rho}{2}\parallel\bm{x}_i-\bm{x}_{j} \parallel ^2).
 \end{aligned}
 \end{equation}
 In \eqref{eq:conventional augemented lagrangian function}, $\bm{x}=[\bm{x}_1^T, \bm{x}_2^T,...,\bm{x}_N^T]^T\in\mathbb{R}^{ND}$ is the augmented states,  $\bm{\lambda}_{i,j}$ is the Lagrange multiplier corresponding to the constraint $\bm{x}_{i}=\bm{x}_{j}$, and all $\bm{\lambda}_{i,j}$ for $e_{i,j}\in E$ are stacked into $\bm{\lambda}\in\mathbb{R}^{|E|D}$. $\rho$ is the penalty parameter, which is a positive constant scalar. 
 
The above ADMM algorithm cannot protect the privacy of participating agents as states are exchanged and disclosed explicitly among neighboring agents.  To facilitate privacy design, we propose a new ADMM with time-varying penalty matrices in the following subsection, which will enable the integration of homomorphic cryptography and decentralized optimization in Sec. \uppercase\expandafter{\romannumeral3}.   
\subsection{ADMM with Time-varying  Penalty Matrices}

Motivated by the fact that  ADMM allows time-varying penalty matrices \cite{kontogiorgis1998variable,he2002new}, we present in the following an ADMM with time-varying penalty matrices.  It is worth noting that \cite{kontogiorgis1998variable,he2002new} deal with a two-block $(N=2)$ problem. While in this paper, we consider a more general problem with $N\ge3$ blocks, whose convergence is more difficult to analyze. The generalization from $N=2$ to $N\ge3$ is highly non-trivial. In fact, as indicated in \cite{chen2016direct}, a direct extension from two-block to multi-block convex minimization is not necessarily convergent. 

We first reformulate \eqref{eq:global function} in a more compact form:
\begin{equation}\label{eq:compact distributed admm form}
	\begin{aligned}
		\mathop {\min }\limits_{\bm{x}} & \qquad f(\bm{x})  \\
		\textrm{subject to} & \qquad A\bm{x}=\bm{0}, 
	\end{aligned}
\end{equation}
where $\bm{x}=[\bm{x}_1^T, \bm{x}_2^T,...,\bm{x}_N^T]^T\in\mathbb{R}^{ND}$, $f(\bm{x})=\sum\limits_{i=1}^{N}f_i(\bm{x}_i)$, and $A=[a_{m,n}]\otimes I_D\in
\mathbb{R}^{|E|D\times ND}$ is the
edge-node incidence matrix of graph $G$ as defined in
\cite{wei2012distributed}, with its $|E|D$
rows corresponding to the $|E|$ communication links and the $ND$ columns corresponding to the $N$ agents. The symbol $\otimes$ denotes Kronecker product. The  $a_{m,n}$ element is defined as
\begin{eqnarray} \nonumber
a_{m,n}=\left\{\begin{matrix}
1& \textrm {if the } m^{th} \textrm{ edge originates from agent }  n,  \\
-1& \textrm {if the } m^{th} \textrm{ edge terminates at agent }  n, \\
0& \textrm{otherwise}.
\end{matrix}\right.
\end{eqnarray}
Here we define that each edge $e_{i,j}$ originates from agent $i$ and
terminates at agent $j$. 

Let $\bm{\lambda}_{i,j}$ be the Lagrange multiplier corresponding to the constraint $\bm{x}_{i}=\bm{x}_{j}$, then we can form an augmented Lagrangian function of problem \eqref{eq:compact distributed admm form} as
\begin{equation} \label{eq:augemented lagrangian function}
\begin{aligned}
\lefteqn{\mathcal{L}(\bm{x},\bm{\lambda}, \bm{\rho})=\sum\limits_{i=1}^{N}f_i(\bm{x}_{i}) }\\
& \qquad+ \sum\limits_{e_{i,j}\in E}(\bm{\lambda}_{i,j}^{T}(\bm{x}_i-\bm{x}_{j})+\frac{\rho_{i,j}}{2}\parallel\bm{x}_i-\bm{x}_{j} \parallel ^2),
\end{aligned}
\end{equation}
 or in a more compact form:
\begin{equation} \label{eq:augemented lagrangian function 2}
\begin{aligned}
\mathcal{L}(\bm{x},\bm{\lambda}, \bm{\rho})=f(\bm{x})+ 
\bm{\lambda}^{T}A\bm{x}+\frac{1}{2}\parallel A\bm{x} \parallel_{\bm{\rho}}^2,
\end{aligned}
\end{equation}
where  
 $\bm{\lambda}=[\bm{\lambda}_{i,j}]_{ij,e_{i,j}\in E}\in\mathbb{R}^{|E|D}$
 is the augmented Lagrange multiplier, $$\bm{\rho}=\mathrm{diag}\{\rho_{i,j}I_D\}_{ij,e_{i,j}\in E}\in\mathbb{R}^{|E|D\times|E|D},\quad \rho_{i,j}>0$$ is the time-varying penalty matrix, and $\parallel A\bm{x} \parallel_{\bm{\rho}}^2=\bm{x}^TA^T\bm{\rho}A\bm{x}.$ 
 
 Note that if $\rho_{i,j}I_D$ is the $m$th block in $\bm{\rho}$, then $e_{i,j}$ is the $m$th edge in $E$, i.e., for the one-dimensional case,  $a_{m,i}=1$ and $a_{m,j}=-1$, and for high dimensional cases,  the $(m,i)$th block of $A$ is $I_D$ and the $(m,j)$th block of $A$ is $-I_D$. 

Now, inspired by \cite{he2002new}, we propose a new ADMM which allows time-varying penalty matrices based on Jacobian update \cite{saad2003iterative}:
   \begin{numcases}{}
 \bm{x}_i^{t+1}=\underset{\bm{x}_i}{\operatorname{argmin}}\mathcal{L}(\bm{x}_1^t,\bm{x}_2^t,...,\bm{x}_{i-1}^t,\bm{x}_i,\bm{x}_{i+1}^t,...,\bm{x}_N^t,\bm{\lambda}^t,\bm{\rho}^t) \nonumber \label{eq:time-varing admm 1}\\
 \qquad\qquad\qquad\quad+\frac{\gamma_i}{2}\parallel \bm{x}_i-\bm{x}_i^t \parallel^2, \\
 \rho_{i,j}^t \to \rho_{i,j}^{t+1}, \label{eq:time-varing admm 2} \\
 \bm{\lambda}_{i,j}^{t+1}=\lambda_{i,j}^{t}+\rho_{i,j}^{t+1}(\bm{x}_i^{t+1}-\bm{x}_j^{t+1}),\quad\forall j\in\mathcal{N}_i  \label{eq:time-varing admm 3}
 \end{numcases}
 for $i=1,2,...,N$. It is worth noting that although the communication graph is undirected, 
 	we introduce both $\bm{\lambda}_{i,j}$ and $\bm{\lambda}_{j,i}$ for $e_{i,j}\in E$ in \eqref{eq:conventional admm 2} and \eqref{eq:time-varing admm 3} to unify the algorithm description.  More specifically, we set $\bm{\lambda}_{i,j}^{0}=\rho_{i,j}^0(\bm{x}_i^{0}-\bm{x}_j^{0})$ at $t=0$ such that $\bm{\lambda}_{i,j}^t=-\bm{\lambda}_{j,i}^t$ holds for all $i=1,2,\cdots, N, j\in\mathcal{N}_i$. In this way, we can unify the update rule of agent $i$ without separating $i>j$ and $i<j$ for $j\in\mathcal{N}_i$, as shown in \eqref{eq:j-ADMM lambda  update}.
\begin{Remark 1}
	 The proximal Jacobian  ADMM \eqref{eq:conventional admm 1}-\eqref{eq:conventional admm 2} can be considered as a special case of \eqref{eq:time-varing admm 1}-\eqref{eq:time-varing admm 3} by assigning the same and constant weight $\rho_{i,j}=\rho$ to different equality constraints $\bm{x}_i=\bm{x}_j$.
	Different from the ADMM which uses the same $\rho$ (which might be time-varying in, e.g., the two-block optimization problem \cite{he1998some}) for all equality constraints, the new approach uses different and time-varying $\rho_{i,j}$ for different equality constraints $\bm{x}_i=\bm{x}_j$. As indicated later, this is key for enabling privacy preservation.  
\end{Remark 1}
\begin{Remark 1}
		We did not use Gauss-Seidel update \cite{wei2012distributed}, which requires a predefined global order and hence as indicated in \cite{deng2017parallel}, is not amenable to parallelism. Different from \cite{deng2017parallel} which has a constant penalty parameter, we intentionally introduce time-varying penalty matrix to enable privacy preservation.  Despite enabling new capabilities  in privacy protection (with the assistance of partially homomorphic Paillier encryption), introducing time-varying penalty matrix also reduces convergence rate to $O(1/t)$, in contrast to the $o(1/t)$  rate in \cite{deng2017parallel}. Besides giving new capabilities in privacy and different result in convergence rate, the novel idea of intentional time-varying penalty matrix also leads to difference in theoretical analysis in comparison with \cite{deng2017parallel}. For example, different from [44] which relies on constant-penalty based monotonically non-increasing sequences to prove convergence, our introduction of time-varying penalty matrix results in non-monotonic sequences which prompted us to use a variational inequality to  facilitate the analysis.
\end{Remark 1}
   It is obvious that the new ADMM \eqref{eq:time-varing admm 1}-\eqref{eq:time-varing admm 3} can be implemented in a decentralized manner. The detailed implementation procedure is outlined in Algorithm \uppercase\expandafter{\romannumeral1}. \\

\begin{flushleft}
	\vspace{-1.3cm}
	\rule{0.49\textwidth}{0.2pt}
\end{flushleft}
\vspace{-0.3cm}

\textbf{Algorithm \uppercase\expandafter{\romannumeral1} }
\begin{flushleft}
	\vspace{-0.95cm}
	\rule{0.49\textwidth}{0.2pt}
\end{flushleft}
\vspace{-0.15cm}

\textbf{Initial Setup:}
Each agent $i$ initializes $\bm{x}_i^{0}$, $\rho_{i,j}^{0}$. 

\textbf{Input:} $\bm{x}_i^{t}$, $\bm{\lambda}_{i,j}^{t-1}$, $\rho_{i,j}^{t}$

\textbf{Output:} $\bm{x}_i^{t+1}$, $\bm{\lambda}_{i,j}^{t}$, $\rho_{i,j}^{t+1}$

\begin{enumerate}
	\item Each agent $i$ sends  $\bm{x}_i^{t}$, $\rho_{i,j}^{t}$ to its neighboring agents, and then set $\rho_{i,j}^t=\min\{\rho_{i,j}^{t},\rho_{j,i}^{t}\}.$ It is clear that $\rho_{i,j}^t=\rho_{j,i}^t$ holds.
	
	\item Each agent $i$ updates $\bm{\lambda}_{i,j}^t$ as follows for $j\in \mathcal{N}_i$
		\begin{eqnarray}
			\bm{\lambda}_{i,j}^{t}=\bm{\lambda}_{i,j}^{t-1}+\rho_{i,j}^t(\bm{x}_{i}^{t}-\bm{x}_{j}^{t}). \label{eq:j-ADMM lambda  update}
		\end{eqnarray}
	It is clear that $\bm{\lambda}_{i,j}^{t}=-\bm{\lambda}_{j,i}^{t}$ holds (note that when $t=0$, we set $	\bm{\lambda}_{i,j}^{0}=\rho_{i,j}^0(\bm{x}_{i}^{0}-\bm{x}_{j}^{0})$). 
	
	\item All agents update  their local vectors in parallel:
	\begin{equation}\label{eq:j-ADMM x update}
		\begin{aligned}
	\lefteqn{	\bm{x}_i^{t+1}\in\textrm{argmin}_{\bm{x}_i}f_i(\bm{x}_{i})+ \frac{\gamma_{i}}{2}\parallel \bm{x}_{i}-\bm{x}_{i}^{t} \parallel ^2}\\
& \qquad+\sum\limits_{j\in\mathcal{N}_i}((\bm{\lambda}_{i,j}^{t})^T\bm{x}_{i}+ \frac{\rho_{i,j}^t}{2}\parallel \bm{x}_{i}-\bm{x}_{j}^{t} \parallel ^2). 
	\end{aligned}
	\end{equation}
	
	Here we added two proximal terms $\frac{\rho_{i,i}^t}{2}\parallel \bm{x}_{i}-\bm{x}_{i}^{t} \parallel ^2$ and $\frac{\gamma_{i}}{2}\parallel \bm{x}_{i}-\bm{x}_{i}^{t} \parallel ^2$ to accommodate the influence of $\bm{x}_i^t$. For all $\gamma_{i}>0$, $\rho_{i,i}^t$ is set to 
	\begin{equation}
	\begin{aligned}
	\rho_{i,i}^t=1-\sum\limits_{j\in\mathcal{N}_i,j\ne i}\rho_{i,j}^t. \label{eq:rhoii}
	\end{aligned}
	\end{equation}
		 
	 \item Each agent $i$ updates $\rho_{i,j}^{t+1}$ for all $j\in\mathcal{N}_i$ and sets $t=t+1$. The detailed update rule for $\rho_{i,j}$ will be elaborated later in Theorem \ref{Theorem:J-ADMM Convergence}.
		
\end{enumerate}

\begin{flushleft}
	\vspace{-1.0cm}
	\rule{0.49\textwidth}{0.2pt}
\end{flushleft}

\begin{Remark 1}
  A weighted ADMM which also assigns different weights to different equality constraints is proposed in  \cite{ling2016weighted}. However, the weights in \cite{ling2016weighted} are constant while Algorithm \uppercase\expandafter{\romannumeral1} allows time-varying  weights in each iteration, which, as shown later, is key to enable the integration of partially homomorphic cryptography with decentralized optimization.
\end{Remark 1}

\subsection{Convergence Analysis}
In this subsection, we rigorously prove the convergence of Algorithm  \uppercase\expandafter{\romannumeral1} under the following standard assumptions \cite{wei2012distributed}:
\begin{Assumption 1} \label{assumption 1}
	Each private local function $f_i: \mathbb{R}^D\to\mathbb{R}$ is convex and continuously differentiable. 
\end{Assumption 1}
\begin{Assumption 1} \label{assumption 2}
	Problem \eqref{eq:compact distributed admm form} has an optimal solution, i.e., the Lagrangian function
	\begin{eqnarray}
	L(\bm{x},\bm{\lambda})=f(\bm{x})+\bm{\lambda}^TA\bm{x} \label{eq:lagrangian function}
	\end{eqnarray} 
   has a saddle point $(\bm{x}^*,\bm{\lambda}^*)$ such that
	$$L(\bm{x}^*,\bm{\lambda})\le L(\bm{x}^*,\bm{\lambda}^*)\le L(\bm{x},\bm{\lambda}^*)$$ holds for all $\bm{x}\in\mathbb{R}^{ND}$ and $\bm{\lambda}\in\mathbb{R}^{|E|D}$. 
\end{Assumption 1}

Denote the iterating results in the $k$th step in Algorithm \uppercase\expandafter{\romannumeral1} as follows:
\begin{displaymath}
\begin{aligned}
&\bm{x}^k=[\bm{x}_1^{kT},\bm{x}_2^{kT},...,\bm{x}_N^{kT} ]^T\in\mathbb{R}^{ND}, \\
&\bm{\lambda}^k=[\bm{\lambda}_{i,j}^k]_{ij,e_{i,j}\in E}\in\mathbb{R}^{|E|D}, \\
&\bm{\rho}^k=\text{diag}\{\rho_{i,j}^kI_D\}_{ij,e_{i,j}\in E}\in\mathbb{R}^{|E|D\times|E|D}. 
\end{aligned}
\end{displaymath}
Further augment the coefficients $\gamma_i$ $(i=1,2,...,N)$ in  \eqref{eq:j-ADMM x update} into the matrix form
\begin{displaymath}
Q_P={\rm diag}\{\gamma_1,\gamma_2,\ldots,\gamma_N\}\otimes I_D\in\mathbb{R}^{ND\times ND},
\end{displaymath}
and augment $\rho_{i,j}^k$ into the following matrix form
\begin{displaymath}
Q_C^k={\rm diag}\{\sum\limits_{j\in\mathcal{N}_1}\rho_{1,j}^k, \sum\limits_{j\in\mathcal{N}_2}\rho_{2,j}^k, \ldots, \sum\limits_{j\in\mathcal{N}_N}\rho_{N,j}^k\}\otimes I_D,
\end{displaymath}
and $Q_C^k\in\mathbb{R}^{ND\times ND}$. By plugging \eqref{eq:rhoii} into $Q_C^k$, we have $Q_C^k=I_{ND}$, i.e., $Q_C^k$  is an identity matrix.

Now we are in position to give the main results of this subsection:
\begin{Theorem 1}\label{Theorem:J-ADMM Convergence}
	Under Assumption \ref{assumption 1} and Assumption  \ref{assumption 2}, Algorithm \uppercase\expandafter{\romannumeral1} is guaranteed to converge to an optimal solution to \eqref{eq:compact distributed admm form} if the following two conditions are met: 
	
	 Condition A:
	 The sequence $\{\bm{\rho}^k\}$ satisfies 
	 $$0\prec\bm{\rho}^0 \preceq\bm{\rho}^k\preceq \bm{\rho}^{k+1}\preceq \bar{\bm{\rho}}, \quad\forall k\ge 0, $$
	 where $\bm{\rho}^0\succ 0$ means that $\bm{\rho}^0$ is positive definite, and similarly $\bm{\rho}^k\preceq \bm{\rho}^{k+1}$ means that $\bm{\rho}^{k+1}-\bm{\rho}^{k}$ is positive semi-definite.
	 
	Condition B:
	$Q_P+Q_C^k\succ A^T\bar{\bm{\rho}}A.$
\end{Theorem 1}	
{\it Proof}: The proof is provided in the Appendix. \hfill $\blacksquare$

\begin{Theorem 1}\label{Theorem: convergence rate}
	The convergence rate of Algorithm \uppercase\expandafter{\romannumeral1} is $O(1/t)$, where $t$ is the iteration time.
\end{Theorem 1}
{\it Proof}: The proof is provided in the Appendix. \hfill $\blacksquare$

\section{Privacy-Preserving Decentralized Optimization}
 Algorithm \uppercase\expandafter{\romannumeral1} requires agents to exchange  and disclose states explicitly in each iteration among neighboring agents to reach  consensus on the final optimal
solution. In this section, we combine partially homomorphic cryptography with  Algorithm \uppercase\expandafter{\romannumeral1} to propose a privacy-preserving approach for decentralized optimization. We first give the definition of privacy used in this paper.

\begin{Definition 1}
	A mechanism $\mathcal{M}: \mathcal{M}(\mathcal{X})\to \mathcal{Y}$ is defined to be privacy preserving if the input  $\mathcal{X}$  cannot be uniquely derived from the output $\mathcal{Y}$.
\end{Definition 1}

This definition of privacy is inspired by the privacy-preservation definitions in  \cite{du2004privacy,liu2006random,han2010privacy,cao2014privacy,lou2017privacy,yan2013distributed} which take advantages of the fact that if a system of equations has infinite number of solutions, it is impossible to derive the exact value of the original input data from the output data. Therefore, privacy preservation is achieved (see, e.g, Part 4.2.2 in \cite{han2010privacy}). Next, we introduce the Paillier cryptosystem and our privacy-preserving approach.

\subsection{Paillier Cryptosystem}
Our method uses the flexibility of time-varying penalty matrices in Algorithm \uppercase\expandafter{\romannumeral1} to enable the incorporation of Paillier cryptosystem \cite{paillier1999public} in a completely decentralized setting.  
 The Paillier cryptosystem is a public-key cryptosystem which uses a pair of keys: a public key and a private key. The public key can be disseminated publicly and used by any person to encrypt a message, but the message can only be decrypted by the private key. The Paillier cryptosystem includes three algorithms, which are detailed below:\\
\begin{flushleft}
	\vspace{-1.3cm}
	\rule{0.49\textwidth}{0.2pt}
\end{flushleft}
\vspace{-0.25cm}
\textbf{Paillier cryptosystem }
\begin{flushleft}
	\vspace{-0.9cm}
	\rule{0.49\textwidth}{0.2pt}
\end{flushleft}
	\vspace{-0.15cm}
\textbf{Key generation:}
\begin{enumerate}
	\item Choose two large prime numbers $p$ and $q$ of equal bit-length and compute $n=pq$.
	\item Let $g=n+1$.
	\item Let $\lambda=\phi(n)=(p-1)(q-1)$, where $\phi(\cdot)$ is the Euler's totient function.
	\item Let $\mu=\phi(n)^{-1}\mod n$ which is the modular multiplicative inverse of $\phi(n)$.
	\item The public key $k_p$ for encryption is $(n,g)$.
	\item  The private key $k_s$ for decryption is $(\lambda,\mu)$.	
\end{enumerate}
\textbf{Encryption ($c=\mathcal{E}(m)$):}

Recall the definitions of $\mathbb{Z}_n=\{z|z\in\mathbb{Z},0\le z<n\}$ and $\mathbb{Z}_n^*=\{z|z\in\mathbb{Z},0\le z<n, \gcd(z,n)=1\}$.
\begin{enumerate}
	\item Choose a random $r\in\mathbb{Z}_n^*$.
	\item The ciphertext is given by $c=g^m\cdot r^n\mod n^2 $, where $m\in\mathbb{Z}_n, c\in\mathbb{Z}_{n^2}^*$. 
\end{enumerate}
\textbf{Decryption ($m=D(c)$):}
\begin{enumerate}
	\item Define the integer division function $L(\mu)=\frac{\mu-1}{n}$.
	\item The plaintext is $m=L(c^\lambda\mod n^2)\cdot \mu \mod n $. 
\end{enumerate}
\begin{flushleft}
	\vspace{-0.85cm}
	\rule{0.49\textwidth}{0.2pt}
\end{flushleft}

A notable feature of Paillier cryptosystem is that it is additively homomorphic, i.e., the ciphertext of $m_1+m_2$ can be obtained from the ciphertext of $m_1$ and $m_2$ directly: 
\begin{eqnarray}
\mathcal{E}(m_1,r_1)\cdot \mathcal{E}(m_2,r_2)=\mathcal{E}(m_1+m_2,r_1r_2), \label{eq:additive}\\
\mathcal{E}(m)^k=\mathcal{E}(km), \quad k\in \mathbb{Z}^+. \qquad
\end{eqnarray} 
Due to the existence of random $r$, the Paillier cryptosystem is resistant to the dictionary attack \cite{goldreich2009foundations}. Since $r_1$ and $r_2$ play no role in the decryption process, \eqref{eq:additive} can be simplified as 
\begin{eqnarray}
\mathcal{E}(m_1)\cdot \mathcal{E}(m_2)=\mathcal{E}(m_1+m_2). 
\end{eqnarray}

\subsection{Privacy-Preserving Decentralized Optimization}
In this subsection, we combine Paillier cryptosystem with Algorithm \uppercase\expandafter{\romannumeral1} to enable privacy preservation in the decentralized solving of optimization problem \eqref{eq:global function}. First, note that solving \eqref{eq:j-ADMM x update} amounts to solving the following problem:
 \begin{equation}\label{eq:j-ADMM x update 1}
 \begin{aligned}
 \triangledown f_i(\bm{x}_{i})
 +\sum\limits_{j\in\mathcal{N}_i}(\bm{\lambda}_{i,j}^{t}+\rho_{i,j}^t(\bm{x}_{i}-\bm{x}_{j}^{t} ))+\gamma_i(\bm{x}_i-\bm{x}_i^t)=\bm{0}. 
 \end{aligned}
 \end{equation}

Let $\bm{\lambda}_{i}=\sum\limits_{j\in N_i}\bm{\lambda}_{i,j}$, then \eqref{eq:j-ADMM x update 1} reduces to the following equation
 \begin{equation}\label{eq:j-ADMM x update 2}
 \begin{aligned}
 \triangledown f_i(\bm{x}_{i}) +(\sum\limits_{j\in\mathcal{N}_i}\rho_{i,j}^t+\gamma_i)\bm{x}_i+\bm{\lambda}_i^t
- \sum\limits_{j\in\mathcal{N}_i}\rho_{i,j}^t\bm{x}_{j}^{t} -\gamma_i\bm{x}_i^t=\bm{0}. 
 \end{aligned}
 \end{equation}
 
Given that we have set $\rho_{i,i}^t=1-\sum\limits_{j\in\mathcal{N}_i,j\ne i}\rho_{i,j}^t$ in \eqref{eq:rhoii}, we can further reduce  \eqref{eq:j-ADMM x update 2} to
  \begin{equation}\label{eq:j-ADMM x update 3}
  \begin{aligned}
\lefteqn{ \triangledown f_i(\bm{x}_{i})+(1+\gamma_i)\bm{x}_i+\bm{\lambda}_i^t}\\
	&\qquad-\sum\limits_{j\in\mathcal{N}_i}\rho_{i,j}^t(\bm{x}_{j}^{t}- \bm{x}_{i}^{t})
 -(1+\gamma_i)\bm{x}_i^t=\bm{0}. 
  \end{aligned}
  \end{equation}

By constructing $\rho_{i,j}^t, i\ne j$ as the product of two random positive numbers, i.e., $\rho_{i,j}^t=b_{i\shortrightarrow  j}^t\times b_{j\shortrightarrow  i}^t=\rho_{j,i}^t$, with $b_{i\shortrightarrow  j}^t$ only known to agent $i$ and  $b_{j\shortrightarrow  i}^t$ only known to agent $j$, we can propose the following privacy-preserving solution to \eqref{eq:global function} based on Algorithm \uppercase\expandafter{\romannumeral1}:
\begin{flushleft}
	\vspace{-0.88cm}
	\rule{0.49\textwidth}{0.2pt}
\end{flushleft}
\vspace{-0.25cm}

\textbf{Algorithm \uppercase\expandafter{\romannumeral2}}
\begin{flushleft}
	\vspace{-0.9cm}
	\rule{0.49\textwidth}{0.2pt}
\end{flushleft}
\vspace{-0.1cm}

\textbf{Initial Setup:}
Each agent initializes $\bm{x}_i^{0}$. 

\textbf{Input:} $\bm{x}_i^{t}$, $\bm{\lambda}_{i,j}^{t-1}$

\textbf{Output:} $\bm{x}_i^{t+1}$, $\bm{\lambda}_{i,j}^{t}$

\begin{enumerate}	
		\item Agent $i$ encrypts $-\bm{x}_i^{t}$ with its public key $k_{pi}$: $$\bm{x}_i^{t}\to\mathcal{E}_i(-\bm{x}_i^{t}).$$ Here the subscript $i$ denotes encryption using the public key of agent $i$.  
		\item Agent $i$ sends $\mathcal{E}_i(-\bm{x}_i^{t})$ and its public key $k_{pi}$ to neighboring agents.
		\item Agent $j\in\mathcal{N}_i$  encrypts $\bm{x}_j^{t}$ with agent $i$'s public key $k_{pi}$:  $$\bm{x}_j^{t}\to\mathcal{E}_i(\bm{x}_j^{t}).$$
		\item Agent $j\in\mathcal{N}_i$ computes the difference directly in ciphertext:
		$$\mathcal{E}_i(\bm{x}_j^{t}-\bm{x}_i^{t})=\mathcal{E}_i(\bm{x}_j^{t})\cdot\mathcal{E}_i(-\bm{x}_i^{t}).$$
		\item Agent $j\in\mathcal{N}_i$ computes the $b_{j\shortrightarrow  i}^t$-weighted difference in ciphertext:
		$$\mathcal{E}_i(b_{j\shortrightarrow  i}^t(\bm{x}_j^{t}-\bm{x}_i^{t}))=(\mathcal{E}_i(\bm{x}_j^{t}-\bm{x}_i^{t}))^{b_{j\shortrightarrow  i}^t}.$$
		\item  Agent $j\in\mathcal{N}_i$ sends $\mathcal{E}_i(b_{j\shortrightarrow  i}^t(\bm{x}_j^{t}-\bm{x}_i^{t}))$ back to agent $i$.
		\item Agent $i$ decrypts the message received from $j$ with its private key $k_{si}$ and multiples the result with $b_{i\shortrightarrow  j}^t$ to get $\rho_{i,j}^t(\bm{x}_j^t-\bm{x}_i^t)$.
		\item  Computing \eqref{eq:j-ADMM lambda  update}, agent $i$ obtains $\bm{\lambda}_{i,j}^{t}$. 
		\item Computing \eqref{eq:j-ADMM x update 3}, agent $i$ obtains $\bm{x}_i^{t+1}$.
	    \item Each agent updates $b_{i\shortrightarrow  j}^t$ to $b_{i\shortrightarrow  j}^{t+1}$ and sets $t=t+1$.
 \end{enumerate}

\begin{flushleft}
	\vspace{-0.85cm}
	\rule{0.49\textwidth}{0.2pt}
\end{flushleft}
Several remarks are in order:
\begin{enumerate}
	\item The only situation that a neighbor knows the state of agent $i$ is when $\bm{x}_i^t=\bm{x}_j^t$ is true for $j\in \mathcal{N}_i$. Otherwise, agent $i$'s state $\bm{x}_i^t$ is encrypted and will not be revealed to its neighbors.
	\item Agent $i$'s state $\bm{x}_i^t$ and its intermediate communication data $b_{j\shortrightarrow  i}^t(\bm{x}_j^t-\bm{x}_i^t)$ will not be revealed to outside eavesdroppers, since they are encrypted.
	\item The state of agent $j\in \mathcal{N}_i$ will not be revealed to agent $i$, because the decrypted message obtained by agent $i$ is $b_{j\shortrightarrow  i}^t(\bm{x}_j^t-\bm{x}_i^t)$ with $b_{j\shortrightarrow  i}^t$ only known to agent $j$ and varying in each iteration.
	\item We encrypt $\mathcal{E}_i(-\bm{x}_i^{t})$ because it is much easier to compute addition in ciphertext. The issue regarding encryption of signed values using Paillier will be addressed in Sec. \uppercase\expandafter{\romannumeral5}.  
	\item Paillier encryption cannot be performed  on vectors directly. For vector messages $\bm{x}_i^t\in\mathbb{R}^D$, each element of the vector (a real number) has to be encrypted separately. For notation convenience, we still denote it in the same way as scalars, e.g., $\mathcal{E}_i(-\bm{x}_i^{t})$.
	\item  Paillier cryptosystem only works for integers, so additional steps have to be taken to convert real values in optimization to integers. This may lead to quantization errors. A common workaround is to scale the real value before quantization, as discussed in detail in Sec. \uppercase\expandafter{\romannumeral5}.
	\item By incorporating Paillier cryptosystem,  it is obvious that the computation complexity and communication load will increase. However, we argue that the privacy provided matters more than this disadvantage when privacy is of primary concern. Furthermore, our experimental results on Raspberry Pi boards confirm that the added communication and computation overhead is fully manageable on embedded microcontrollers (cf. Sec.  \uppercase\expandafter{\romannumeral7}). 
	 \item Our approach is more suitable for small and medium sized optimization problems such as the source localization problem \cite{zhang2017distributed}  and  power system monitoring  problem \cite{wang2018distributed} addressed in our prior work. 
\end{enumerate}

The key to achieve privacy preservation is to construct $\rho_{i,j}^t, i\ne j$ as the product of two random positive numbers $b_{i\shortrightarrow  j}^t$ and $b_{j\shortrightarrow  i}^t$, with $b_{i\shortrightarrow  j}^t$ generated by and only known to agent $i$ and  $b_{j\shortrightarrow  i}^t$ generated by and only known to agent $j$. Next we show that the privacy preservation mechanism does not affect the convergence to the optimal solution.
\begin{Theorem 1} \label{Theorem: alggorithm 2 convergence}
The privacy-preserving algorithm \uppercase\expandafter{\romannumeral2} will generate a solution in an $\varepsilon$ ball around the optimum if $b_{i\shortrightarrow  j}^t$,  $b_{j\shortrightarrow  i}^t$, and $\gamma_i$ are updated in the following way (where $\varepsilon$ depends on the quantization error):
	\begin{enumerate}
		\item  $b_{i\shortrightarrow  j}^t$ is randomly chosen from $[b_{i\shortrightarrow  j}^{t-1},\bar{b}_{i\shortrightarrow  j}]$, with $\bar{b}_{i\shortrightarrow  j}>0$ denoting a predetermined constant only known to agent $i$;
		\item $\gamma_i$ is chosen randomly in the interval $[N\bar{b}^2, \bar{\bar{b}}]$, with  $\bar{b}>\max\{\bar{b}_{i\shortrightarrow  j}\}$ denoting a predetermined positive constant known to everyone and $\bar{\bar{b}}$ a threshold chosen arbitrarily by agent $i$ and only known to agent $i$. 
	\end{enumerate}
\end{Theorem 1}	
{\it Proof:} It can be easily obtained that if $b_{i \shortrightarrow j}^t$ is updated following 1), and $\gamma_i$ is updated following 2), then Condition A and Condition B in Theorem \ref{Theorem:J-ADMM Convergence} will be met automatically. Therefore, the states in algorithm II should converge to the optimal solution. However, since Paillier cryptosystem only works on unsigned integers, it requires converting real-valued states  to integers using e.g., fixed-point arithmetic encoding \cite{pythonpaillier} (after scaled by a large number $N_{\max}$, cf. Sec. V), which leads to quantization errors. The quantization errors lead to numerical errors on the final solution and hence the ``$\varepsilon$-ball" statement in Theorem \ref{Theorem: alggorithm 2 convergence}.  It is worth noting that the numerical error here is no different from the conventional quantization errors met by all algorithms when implemented in practice on a computer. A quantized analysis of the $\varepsilon$-ball is usually notoriously involved and hence we refer interested readers to \cite{zhu2016quantized} which is dedicated to this problem. Furthermore, we would like to emphasize that this quantization error can be made arbitrarily small by using an arbitrarily large $N_{\max}$. In fact, our simulation results in Sec. VI B showed  that under $N_{\max}=10^6$, the final error was on the order of $10^{-14}$.  \hfill $\blacksquare$

\section{Privacy Analysis}
As indicated in the introduction, our approach aims to protect the privacy of agents' intermediate  states $\bm{x}_i^t$s and gradients of $f_i$s as well as the objective functions.
In this section, we rigorously prove that these private information cannot
be inferred by honest-but-curious adversaries and external eavesdroppers, which
are commonly used attack models in privacy studies \cite{li2010secure}
(cf. definition in Sec. \uppercase\expandafter{\romannumeral1}).
 It is worth noting that the form of each agent's local objective function can also be totally blind to others, e.g., whether it is a quadratic, exponential, or other forms of convex functions is only known to an agent itself. 

As indicated in Sec. \uppercase\expandafter{\romannumeral3},  our approach in Algorithm \uppercase\expandafter{\romannumeral2} guarantees that state information is not leaked to any neighbor in one iteration. However, would some information get leaked over time? More specifically, if an honest-but-curious adversary observes carefully its communications with  neighbors over several steps, can it put together all the received information to infer its neighbor's state?

We can rigorously prove that an honest-but-curious adversary cannot infer the exact states of its neighbors even by collecting samples from multiple steps. 

\begin{Theorem 1} \label{theorem:privacy}
	Assume that all agents follow Algorithm \uppercase\expandafter{\romannumeral2}.  Then agent $j$'s exact state value $\bm{x}_j^k$ cannot be inferred by an honest-but-curious  agent $i$ unless $\bm{x}_i^k=\bm{x}_j^k$ is true.  
\end{Theorem 1}

{\it Proof:} Suppose that an honest-but-curious agent $i$ collects information from $K$ iterations to infer the information of a neighboring agent $j$.  From the perspective of adversary agent $i$, the measurements (corresponding to neighboring agent $j$) seen in each iteration $k$ are $\bm{y}^k=b_{i\shortrightarrow  j}^kb_{j\shortrightarrow  i}^k(\bm{x}_j^k-\bm{x}_i^k)$ $(k=0,1,...,K)$, i.e., adversary agent $i$ can establish $(K+1)D$ equations based on received information:
 \begin{eqnarray}
 \left\{
 \begin{aligned}
 \bm{y}^0&=b_{i\shortrightarrow  j}^0b_{j\shortrightarrow  i}^0(\bm{x}_j^0-\bm{x}_i^0), \\
 \bm{y}^1&=b_{i\shortrightarrow  j}^1b_{j\shortrightarrow  i}^1(\bm{x}_j^1-\bm{x}_i^1), \\
 &\quad\vdots\label{eq:equation set}  \\
 \bm{y}^{K-1} &=b_{i\shortrightarrow  j}^{K-1}b_{j\shortrightarrow  i}^{K-1}(\bm{x}_j^{K-1}-\bm{x}_i^{K-1}), \\
 \bm{y}^K&=b_{i\shortrightarrow  j}^Kb_{j\shortrightarrow  i}^{K}(\bm{x}_j^K-\bm{x}_i^K). 
 \end{aligned}
 \right.
 \end{eqnarray}
 
To the adversary agent $i$, in the system of equations \eqref{eq:equation set}, $\bm{y}^k, b_{i\shortrightarrow  j}^k, \bm{x}_i^k$ $(k=0,1,2,...,K)$ are known, but $\bm{x}_j^k, b_{j\shortrightarrow  i}^k$ $(k=0,1,2,...,K)$ are unknown. So the above system of  $(K+1)D$ equations contains $(K+1)D+K+1$ unknown variables. It is clear that adversary agent $i$ cannot
solve the system of equations \eqref{eq:equation set} to infer the exact values of unknowns $\bm{x}_j^k$ and $b_{j\shortrightarrow  i}^k$ $(k=0,1,2,...,K)$ of agent $j$. It is worth noting that if for some time index $k$, $\bm{x}_j^k =\bm{x}_i^k$ happens to be true, then adversary agent $i$ will be able to know that agent $j$ has the same state at this time index based on the fact that $\bm{y}^k$ is $\bm{0}$. \hfill $\blacksquare$

Using a similar way of reasoning, we can obtain that an honest-but-curious adversary agent $i$ cannot  infer the exact gradient of objective function $f_j$ from a neighboring agent $j$ at any point when agent $j$ has another legitimate neighbor other than the honest-but curious neighbor $i$.
 
 \begin{Theorem 1} \label{corollary:gradient privacy}
 In Algorithm \uppercase\expandafter{\romannumeral2}, the exact gradient of $f_j$ at any point cannot be inferred by an honest-but-curious agent $i$ if agent $j$ has another legitimate neighbor.
 \end{Theorem 1}
  {\it Proof:} Suppose that an honest-but-curious adversary agent $i$ collects information from $K$ iterations to infer the gradient of function $f_j$ of a neighboring agent $j$. The adversary agent $i$ can establish $KD$ equations corresponding to  the gradient of $f_j$ by making use of the fact that the update rule \eqref{eq:j-ADMM x update 3} is publicly known, i.e.,
 \begin{equation} \label{eq:j} 
 \left\{
 \begin{aligned}
 \lefteqn{ \triangledown f_j(\bm{x}_{j}^{1})+(1+\gamma_j)\bm{x}_j^{1}+\bm{\lambda}_j^0}\\
 &\qquad\qquad\qquad-\sum\limits_{m\in\mathcal{N}_j}\rho_{j,m}^0(\bm{x}_{m}^{0}- \bm{x}_{j}^{0})
 -(1+\gamma_j)\bm{x}_j^0=\bm{0},  \\
 \lefteqn{ \triangledown f_j(\bm{x}_{j}^{2})+(1+\gamma_j)\bm{x}_j^{2}+\bm{\lambda}_j^1}\\
 &\qquad\qquad\qquad-\sum\limits_{m\in\mathcal{N}_j}\rho_{j,m}^1(\bm{x}_{m}^{1}- \bm{x}_{j}^{1})
 -(1+\gamma_j)\bm{x}_j^1=\bm{0},  \\
 &\qquad\qquad\qquad\qquad\vdots \\
 \lefteqn{ \triangledown f_j(\bm{x}_{j}^{K-1})+(1+\gamma_j)\bm{x}_j^{K-1}+\bm{\lambda}_j^{K-2}}\\
 &\qquad-\sum\limits_{m\in\mathcal{N}_j}\rho_{j,m}^{K-2}(\bm{x}_{m}^{K-2}- \bm{x}_{j}^{K-2})
 -(1+\gamma_j)\bm{x}_j^{K-2}=\bm{0}, \\
 \lefteqn{ \triangledown f_j(\bm{x}_{j}^{K})+(1+\gamma_j)\bm{x}_j^{K}+\bm{\lambda}_j^{K-1}}\\
 &\qquad-\sum\limits_{m\in\mathcal{N}_j}\rho_{j,m}^{K-1}(\bm{x}_{m}^{K-1}- \bm{x}_{j}^{K-1})
 -(1+\gamma_j)\bm{x}_j^{K-1}=\bm{0}.  \\
 \end{aligned}
 \right.
 \end{equation}
  In the system of $KD$ equations \eqref{eq:j}, $\triangledown f_j(\bm{x}_j^k)$ $(k=1,2,...,K)$, $\gamma_j$, and $\bm{x}_{j}^{k}$ $(k=0,1,2,...,K)$ are unknown to adversary agent $i$. 
  Parameters $\bm{\lambda}_j^k$ and $\sum\limits_{m\in\mathcal{N}_j}\rho_{j,m}^k(\bm{x}_{m}^{k}- \bm{x}_{j}^{k})$   $(k=0,1,2,...,K-1)$ are known to adversary agent $i$ only when agent $j$ has agent $i$ as the only neighbor. Otherwise, $\bm{\lambda}_j^k$ and $\sum\limits_{m\in\mathcal{N}_j}\rho_{j,m}^k(\bm{x}_{m}^{k}- \bm{x}_{j}^{k})$ $(k=0,1,2,...,K-1)$ are unknown to adversary agent $i$. Noting that $\bm{\lambda}_j^{k+1}=\bm{\lambda}_j^{k}-\sum\limits_{m\in\mathcal{N}_j}\rho_{j,m}^{k+1}(\bm{x}_{m}^{k+1}- \bm{x}_{j}^{k+1})$ and  $\bm{\lambda}_j^0=-\sum\limits_{m\in\mathcal{N}_j}\rho_{j,m}^0(\bm{x}_{m}^{0}- \bm{x}_{j}^{0})$, we can see that the above system of $KD$ equations contains 
  $3KD+D+1$ unknowns when agent $j$ has more than one neighbor. Therefore, adversary agent $i$ cannot infer the exact values of $\triangledown f_j(\bm{x}_{j}^{k})$ by solving \eqref{eq:j}.  
   
  It is worth noting that after the optimization converges, adversary agent $i$ can have another piece of information according to the KKT conditions \cite{deng2017parallel}:
   \begin{equation} \label{eq:optimal}
   \begin{aligned}
   \triangledown f_j(\bm{x}_{j}^{*})=-\bm{\lambda}_j^*
   \end{aligned}
   \end{equation}
   where $\bm{x}_{j}^{*}$ denotes the optimal solution and $\bm{\lambda}_j^*$ denotes the optimal multiplier. However, since  $\bm{\lambda}_j^*$ is known to adversary agent $i$ only when agent $j$ has agent $i$ as the only neighbor, we have that adversary agent $i$ cannot infer the exact value of $f_j$ at any point when agent $j$ has another legitimate neighbor besides an honest-but curious neighbor $i$.  \hfill $\blacksquare$
   
   	Using a similar way of reasoning, we have the following corollary corresponding to the situation where agent $j$ has honest-but-curious agent $i$ as the only neighbor.
  
   \begin{Corollary 1} \label{corollary:gradient privacy one neighbor}
   	 In Algorithm \uppercase\expandafter{\romannumeral2}, the exact gradient of $f_j$ at the optimal solution can be inferred by an honest-but-curious agent $i$ if agent $j$ has adversary agent $i$ as the only neighbor. However, at any other point, the gradient of $f_j$ is uninferrable by the adversary agent $i$.
   \end{Corollary 1}
    	 {\it Proof:} Following a similar line of reasoning of Theorem \ref{corollary:gradient privacy}, we can obtain the above Corollary. \hfill $\blacksquare$
    
 Based on Theorem \ref{theorem:privacy}, Theorem \ref{corollary:gradient privacy}, and Corollary \ref{corollary:gradient privacy one neighbor}, we can obtain that agent $i$ cannot infer agent $j$'s local objective function $f_j$.
 \begin{Corollary 1}	
  In Algorithm \uppercase\expandafter{\romannumeral2}, agent $j$'s local objective function $f_j$ cannot be inferred by an honest-but-curious agent $i$.
  \end{Corollary 1}
   {\it Proof:} According to Theorem \ref{theorem:privacy}, Theorem \ref{corollary:gradient privacy}, and  Corollary \ref{corollary:gradient privacy one neighbor}, the intermediate states and corresponding gradients of $f_j$ cannot be inferred by adversary $i$. Therefore, adversary $i$ cannot infer agent $j$'s local objective function $f_j$ as well.  \hfill $\blacksquare$

 Furthermore, we have that an external eavesdropper cannot infer any private information of all agents. 
\begin{Corollary 1}	
 		All agents' intermediate states, gradients of objective functions, and objective functions cannot be inferred by an external eavesdropper.
 	\end{Corollary 1}
 	{\it Proof:} Since all exchanged messages are encrypted and that cracking the encryption is practically infeasible \cite{goldreich2009foundations}, an external eavesdropper cannot learn anything by intercepting these messages. Therefore, it cannot infer any agent's intermediate states, gradients of objective functions, and objective functions.  \hfill $\blacksquare$

From the above analysis, it is obvious that agent $j$'s private information cannot be uniquely derived by adversaries. However, an honest-but-curious neighbor $i$ can still get some range information about the state $\bm{x}_j^k$ and this range information will become tighter as $\bm{x}_j^k$ converges to the optimal solution as $k\to\infty$ (cf. the simulation results in Fig. \ref{fig:x estimate}). We argue that this is completely unavoidable for any privacy-preserving approaches because all agents have to agree on the same final state, upon which the privacy of $\bm{x}_j^k$ will disappear. In fact, this is  also acknowledged in \cite{huang2015differentially}, which shows that the privacy of $\bm{x}_j^k$ will vanish as $k\to\infty$ and the noise variance converges to zero at the state corresponding to the optimal solution. 
It is worth noting that when the constraint is of a form different from consensus, it may be possible to protect the privacy of $\bm{x}_j^k$ when $k\to \infty$. However, how to incorporate the proposed privacy mechanism in decentralized optimization under non-consensus constraint is difficult and will be addressed in our future work.

 \begin{Remark 1}
   It is worth noting that an adversary agent $i$ can combine systems of equations \eqref{eq:equation set} and \eqref{eq:j}  to infer the information of a neighboring agent $j$.   However, this will not increase the ability of adversary agent $i$  because the combination will not change the fact that the number of unknowns is greater than the number of establishable relevant equations.  In addition, if all other agents collude to infer $\bm{x}_j^k$ of agent $j$, these agents can be considered as one agent which amounts to having a network consisting of two agents. 
 \end{Remark 1}
 
\begin{Remark 1}
	From Theorem \ref{theorem:privacy}, we can see that in decentralized optimization, an agent's information will not be disclosed to other agents no matter how many neighbors it has. This is in distinct difference from the average consensus problem in \cite{ruan2017secure,mo2017privacy} where privacy cannot be protected for an agent if it has an honest-but-curious adversary as the only neighbor.
\end{Remark 1}

\section{Implementation Details}
In this section, we discuss several technical issues that have to be addressed in the implementation of  Algorithm \uppercase\expandafter{\romannumeral2}.
\begin{enumerate}
	\item In modern communication, a real number is represented by a floating point number, while encryption techniques only work for unsigned integers. To deal with this problem, we uniformly multiplied each element of the vector message $\bm{x}_i^t\in\mathbb{R}^D$ (in floating point representation) by a sufficiently large number $N_{\max}$ and round off the fractional part during the encryption to convert it to an integer. After decryption, the result is divided by $N_{\max}$. This process is conducted in each iteration and this quantization brings an error upper-bounded by $\frac{1}{N_{\max}}$. In implementation, $N_{\max}$ can be chosen according to the used data structure. 
	\item As indicated in 1), encryption techniques only work for unsigned integers. In our implementation all integer values are stored in fix-length integers (i.e., long int in	C) and negative values are left in 2's complement format. Encryption and intermediate computations are carried out as if the underlying data were unsigned. When the final message is decrypted, the overflown bits (bits outside the fixed length)	are discarded and the remaining binary number is treated as a signed integer which is later converted back to a real value.
\end{enumerate}

\section{Numerical Experiments}
In this section, we first illustrate the efficiency of the proposed approach using C/C++ implementations. Then we compare our approach with the algorithm in \cite{huang2015differentially} and the algorithm in \cite{lou2017privacy}. The open-source C implementation of the Paillier cryptosystem \cite{bethencourtadvanced} is used in our simulations. We conducted numerical experiments on the following global objective function
\begin{eqnarray}
f(\bm{x})=\sum\limits_{i=1}^{N}\frac{1}{p_i}\parallel H_i\bm{x}-\bm{\theta}_i\parallel^2, \label{eq:simulation}
\end{eqnarray}
which makes the optimization problem \eqref{eq:global function} become
\begin{eqnarray}
\min\limits_{\bm{x}}\qquad \sum\limits_{i=1}^{N}\frac{1}{p_i}\parallel H_i\bm{x}-\bm{\theta}_i\parallel^2 \label{eq:simulation 1}
\end{eqnarray}
with $\bm{\theta}_i\in\mathbb{R}^D$, $H_i=h_i\bm{I}_D$ ($h_i\in\mathbb{R}$), and $p_i>0$ ($p_i\in\mathbb{R}$). Hence, each agent $i$  deals with a private local objective function
\begin{eqnarray} \label{eq:individual function}
f_i(\bm{x}_i)=\frac{1}{p_i}\parallel H_i\bm{x}_i-\bm{\theta}_i\parallel^2, \forall i\in\{1,2,\ldots,N\}. 
\end{eqnarray}
We used the above function \eqref{eq:simulation} because it is easy to verify whether the obtained solution  is the minimal value of the original optimization problem, which should be $\frac{\sum_{i=1}^{N}\frac{2h_i}{p_i}\bm{\theta}_i}{\sum_{i=1}^{N}\frac{2h_i^2}{p_i}}$. Furthermore, \eqref{eq:simulation} makes it easy to compare with \cite{huang2015differentially}, whose verification is also based on \eqref{eq:simulation}.  

In the implementation, the parameters are set as follows: $N_{\max}$ was set to $10^6$ to convert each element in $\bm{x}_i$ to a 64-bit integer during intermediate computations. $b_{i\shortrightarrow  j}^t$ was also scaled up in the same way and represented by a 64-bit integer. The encryption and decryption keys were chosen as 256-bit long.
\subsection{Evaluation of Our Approach}
\begin{figure}[!htbp]
	\begin{center}
		\includegraphics[width=0.42\textwidth]{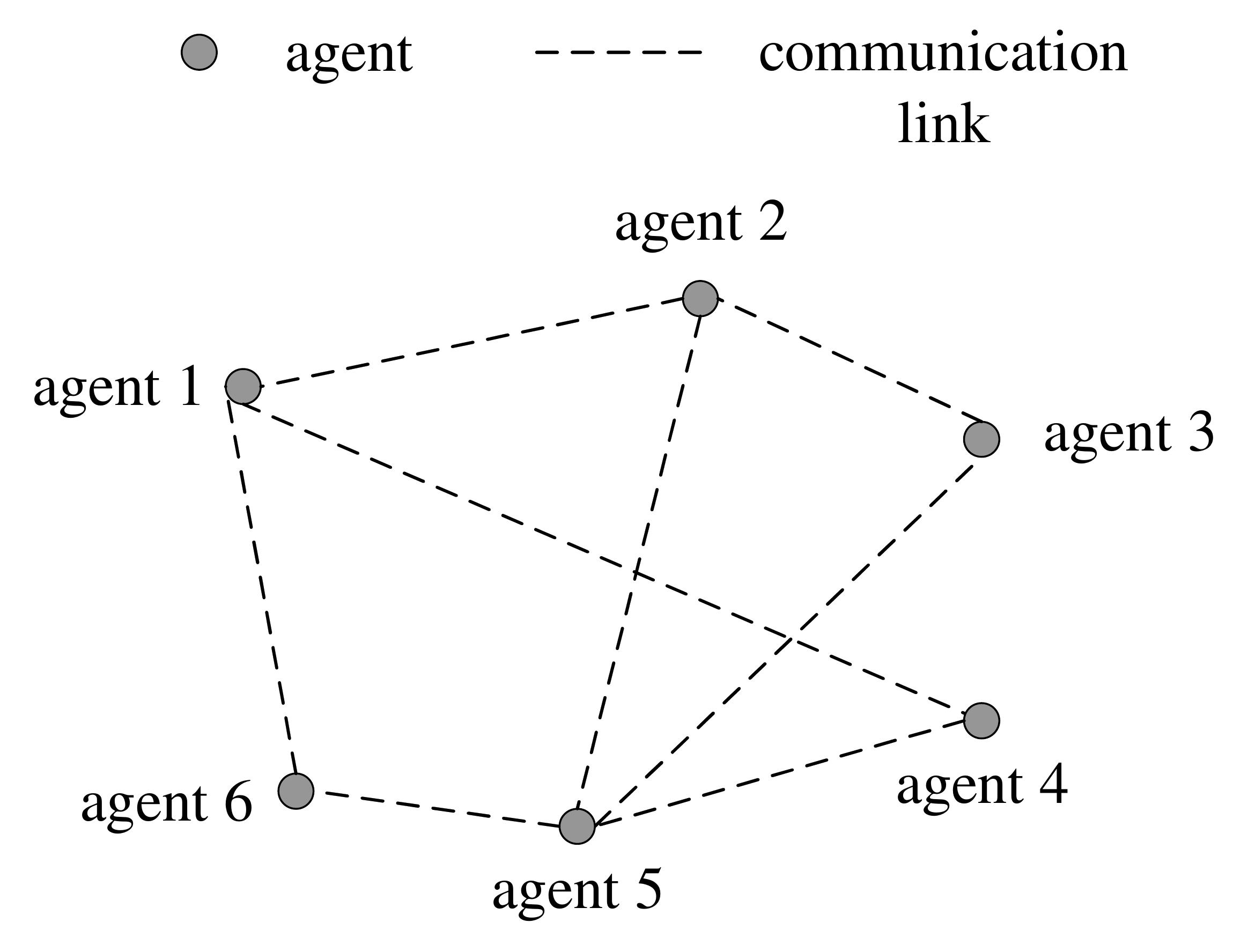}
	\end{center}
	\caption{A network of six agents ($N=6$).}
	\label{fig:network structure illustration}
\end{figure}

\begin{figure}[!htbp]
	\begin{center}
		\includegraphics[width=0.49\textwidth]{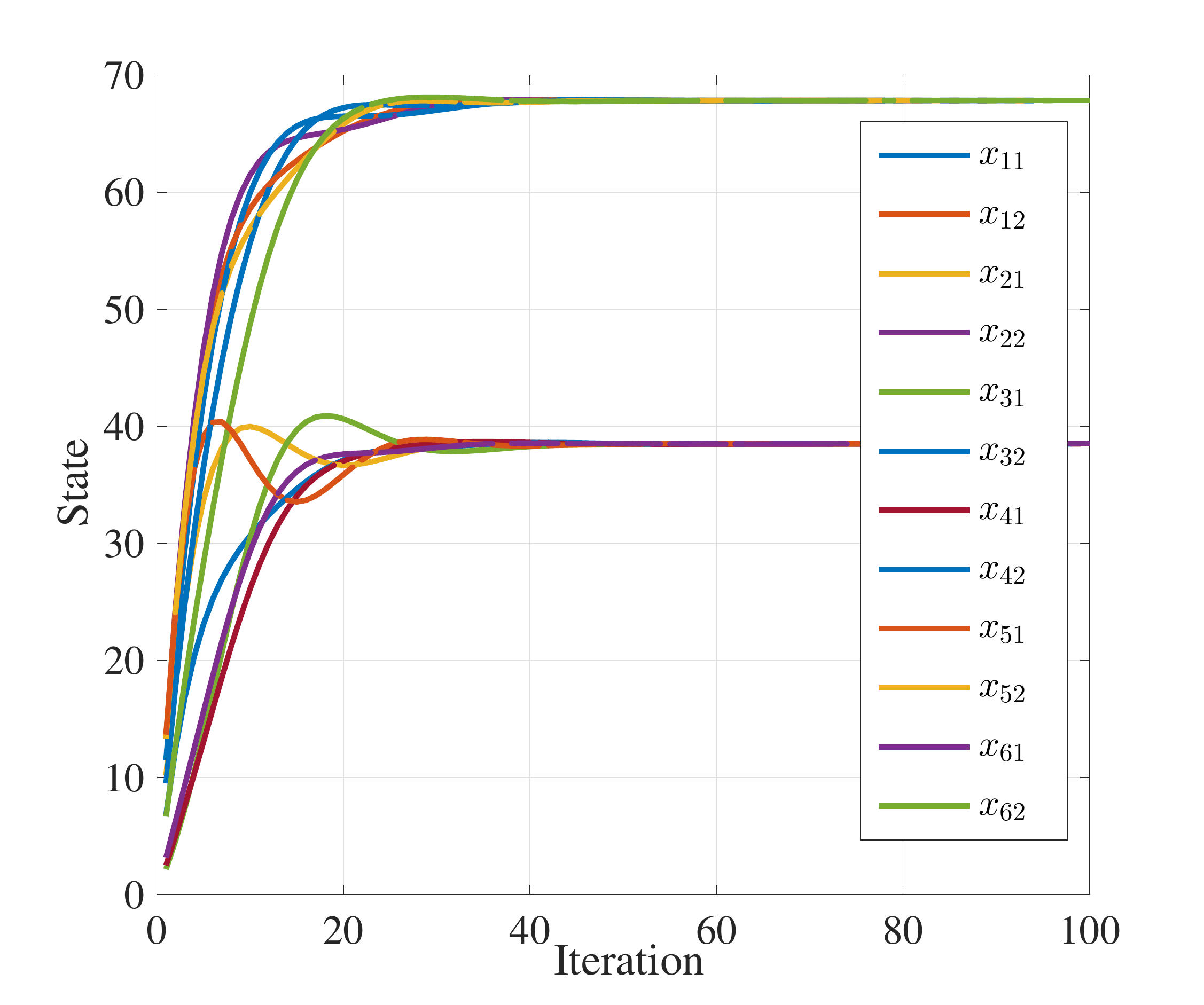}
	\end{center}
	\caption{The evolution of $\bm{x}_i$ ($i=1,2,...,6$).}
	\label{fig:x}
\end{figure}

We implemented Algorithm \uppercase\expandafter{\romannumeral2} on different network topologies, all of which gave the right optimal solution.  Simulation results confirmed that our approach always converged to the optimal solution of  \eqref{eq:simulation 1}. Fig. \ref{fig:x} visualizes the evolution of $\bm{x}_i$ $(i=1,2,...,6)$ in one specific run where the network deployment is illustrated in Fig. \ref{fig:network structure illustration}. In Fig. \ref{fig:x}, $x_{ij}$ $(i=1,2,...,6, j=1,2)$ denotes the $j$th element of $\bm{x}_i$. All $\bm{x}_i$ $(i=1,2,...,6)$ converged to the optimal solution $[38.5;\frac{407}{6}]$. In this run, $\bar{b}$ was set to 0.65 and $\gamma_i$s were set to 3. 

\begin{figure}[!htbp]
	\begin{center}
		\includegraphics[width=0.49\textwidth]{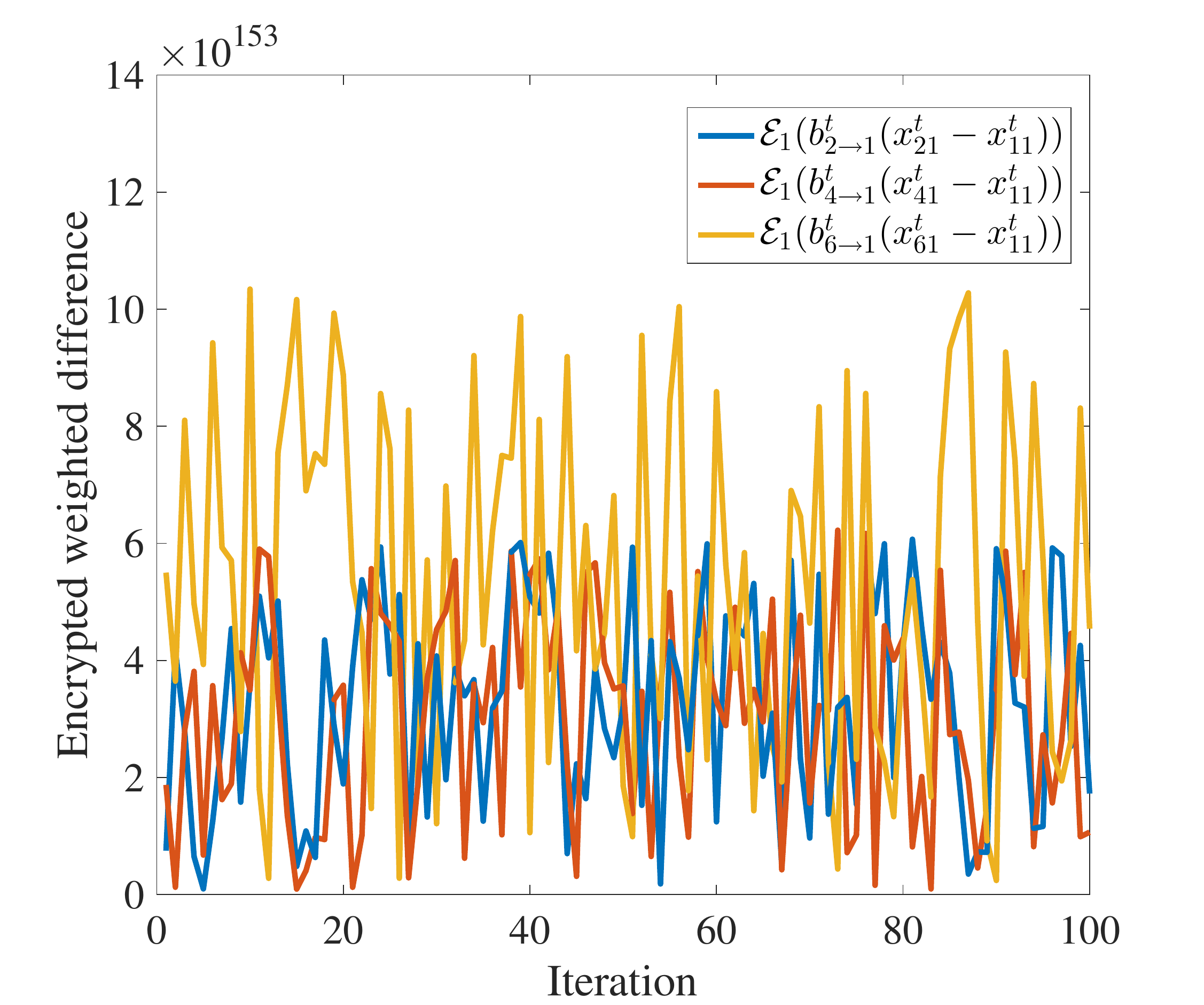}
	\end{center}
	\caption{The evolution of the encrypted weighted differences (in ciphertext)  $\mathcal{E}_1(b_{2 \shortrightarrow  1}^t(\bm{x}_{21}^{t}-\bm{x}_{11}^{t}))$, $\mathcal{E}_1(b_{4 \shortrightarrow  1}^t(\bm{x}_{41}^{t}-\bm{x}_{11}^{t}))$, and $\mathcal{E}_1(b_{6 \shortrightarrow  1}^t(\bm{x}_{61}^{t}-\bm{x}_{11}^{t}))$.}
	\label{fig:ciphertext}
\end{figure}

Fig. \ref{fig:ciphertext} visualizes the encrypted weighted differences (in ciphertext) $\mathcal{E}_1(b_{2\shortrightarrow  1}^t(\bm{x}_{2 1}^{t}-\bm{x}_{11}^{t}))$, $\mathcal{E}_1(b_{4 \shortrightarrow  1}^t(\bm{x}_{41}^{t}-\bm{x}_{11}^{t}))$, and $\mathcal{E}_1(b_{6\shortrightarrow  1}^t(\bm{x}_{61}^{t}-\bm{x}_{11}^{t}))$. It is worth noting that although the states of all agents have converged after about 40 iterations, the encrypted weighted differences (in ciphertext) still appeared random to an outside eavesdropper.

 We also simulated an honest-but-curious adversary who tries to estimate its neighbors' intermediate states and gradients in order to estimate the objective function. We considered the worse case of two agents (A and B) where agent B is the honest-but-curious adversary and intends to estimate the objective function $f_A$ of agent A.  
	  The individual local objective functions are the same as \eqref{eq:individual function} with $\theta_i\in\mathbb{R}$.   Because agent B knows the constraints on agent A's generation of $b_{A\shortrightarrow B}^t$ and $\gamma_A$ (cf. Theorem \ref{Theorem: alggorithm 2 convergence}), it generates estimates of  $b_{A\shortrightarrow B}^t$ and $\gamma_A$ in the same random way.  Then it obtained a series of estimated $x_A^{t}$ and $\triangledown f_A(x_{A}^{t})$ according to \eqref{eq:j}. Finally, agent B used the estimated $x_A^{t}$ and $\triangledown f_A(x_{A}^{t})$ to estimate $f_A$.
 
\begin{figure}[!htbp]
	\begin{center}
		\includegraphics[width=0.49\textwidth]{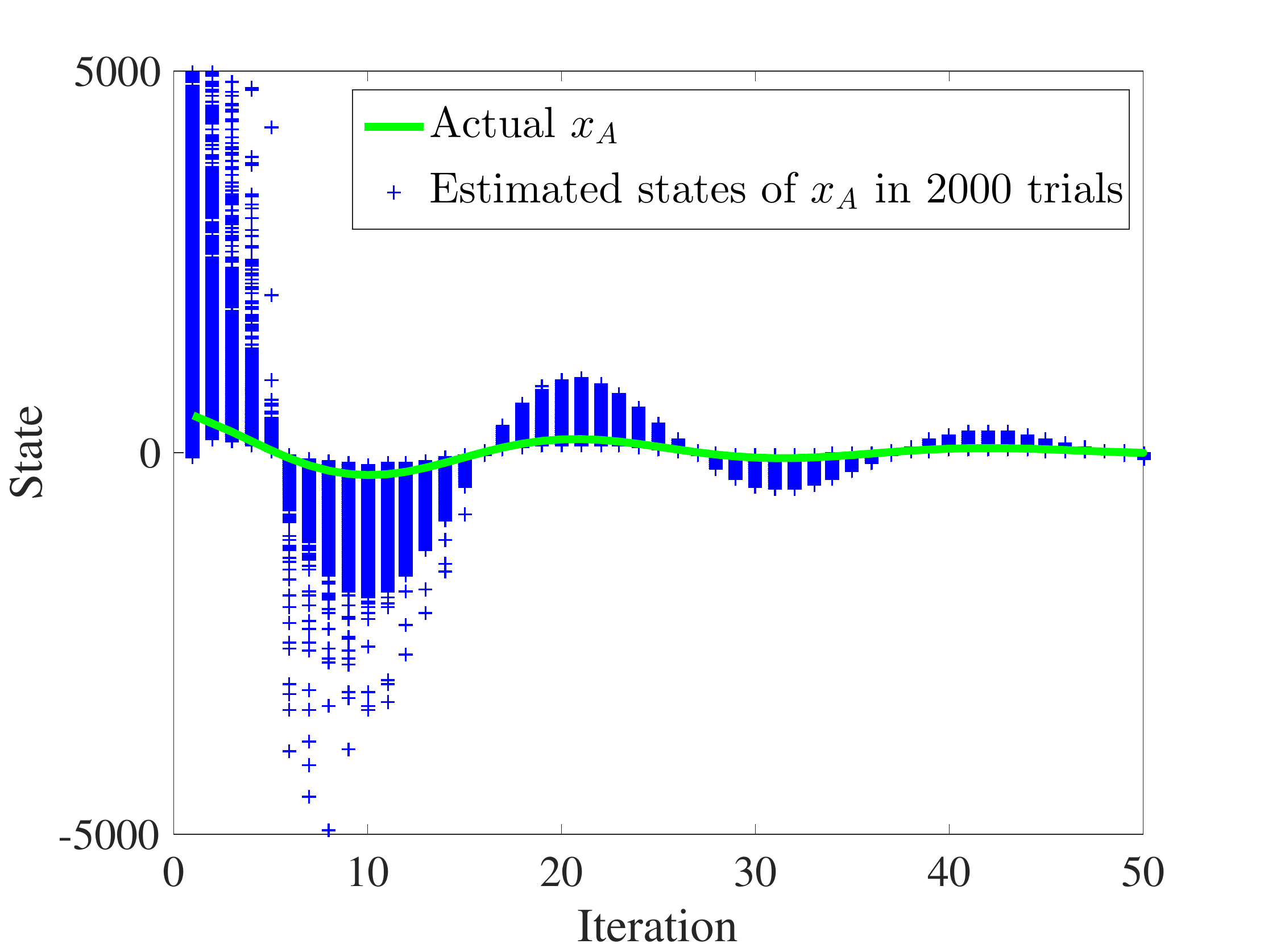}
	\end{center}
	\caption{Estimated states of $x_A$ in 2,000 trials.}
	\label{fig:x estimate}
\end{figure}

\begin{figure}[!htbp]
 	\begin{center}
 		\includegraphics[width=0.49\textwidth]{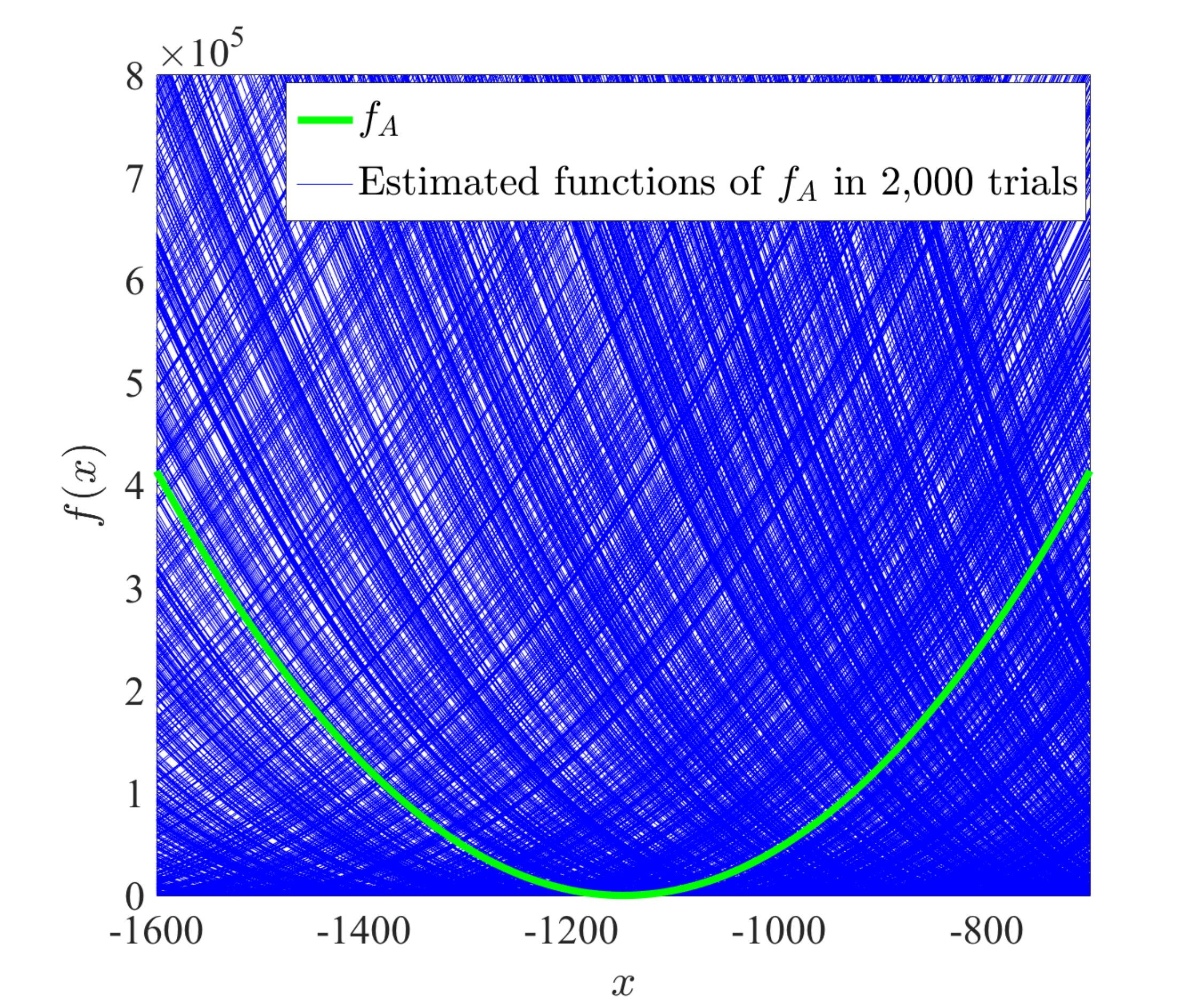}
 	\end{center}
 	\caption{Estimated functions of $f_A$ in 2,000 trials.}
 	\label{fig:estimate}
 \end{figure}

Fig. \ref{fig:x estimate} and Fig. \ref{fig:estimate} show the estimated $x_A$ and $f_A$ in 2,000 trials  when agent B used simple linear regression to estimate $\triangledown f_A(x)$. Fig. \ref{fig:estimate} suggests that agent B cannot get a good estimate of $f_A$. Moreover, it is worth noting that all these estimated functions give the same optimal solution as $f_A$ to the optimization problem \eqref{eq:simulation 1}. 

In addition, the encryption/decryption computation took about 1ms for each agent to communicate with one neighbor at each iteration on a 3.6 GHz CPU, which is manageable in small or medium sized real-time optimization problems such as the source localization problem  \cite{zhang2017distributed} and  power system monitoring problem  \cite{wang2018distributed}  addressed in our prior work. For large sized optimization problems like machine learning with extremely large dimensions, the approach may be computationally too heavy due to the underlying Paillier encryption scheme.

\subsection{Comparison with the algorithm in \cite{huang2015differentially} }
We then compared our approach with the differential-privacy based privacy-preserving optimization algorithm in \cite{huang2015differentially}. Under the communication topology in Fig. \ref{fig:network structure illustration}, we simulated the algorithm in \cite{huang2015differentially} under seven different privacy levels: $\epsilon=0.2, 1, 10, 20, 30, 50, 100$. The global function we used for comparison was  \eqref{eq:simulation} with $p_i$ $(i=1,2,..,6)$ fixed to $2$, $h_i$ $(i=1,2,..,6)$ fixed to $1$, and $\bm{\theta}_i=[0.1\times (i-1)+0.1;0.1\times (i-1)+0.2]$. The domain of optimization was set to $\mathcal{X}=\{(x,y)\in\mathbb{R}^2|x^2+y^2\le 1\}$ for the algorithm in \cite{huang2015differentially}. Note that the optimal solution $[0.35;0.45]$ resided in $\mathcal{X}$. 
    Parameter settings for the algorithm in \cite{huang2015differentially} are detailed as follows: $n=2$, $c=0.5$, $q=0.8$, $p=0.9$, and 
    \begin{eqnarray} \label{eq: aij}
    a_{ij}=\left\{
    \begin{aligned}
    0.2 \qquad\qquad\qquad &  j\in\mathcal{N}_i\verb|\|i, \\
    0 \qquad\qquad\qquad &  j\notin\mathcal{N}_i, \\
    1-\sum\limits_{j\in\mathcal{N}_i\verb|\|i}a_{ij} \qquad & i=j,
    \end{aligned}
    \right.
    \end{eqnarray} 
    for $i=1,2,...,6$. Here $\mathcal{N}_i\verb|\|i$ denotes all values except $i$ in set $\mathcal{N}_i$. Furthermore, we used the performance index $d$ in \cite{huang2015differentially} to quantify the optimization error, which was computed as the average value of squared distances with respect to the optimal solution over $M$ runs \cite{huang2015differentially}, i.e.,
    $$d=\frac{\sum\limits_{i=1}^{6}\sum\limits_{k=1}^{M}\parallel \bm{x}_i^k-[0.35;0.45]\parallel^2}{6M}$$
    with $\bm{x}_i^k$ the obtained solution of agent $i$ in the $k$th run. 
    
\begin{figure}[!htbp]
    	\begin{center}
    		\includegraphics[width=0.49\textwidth]{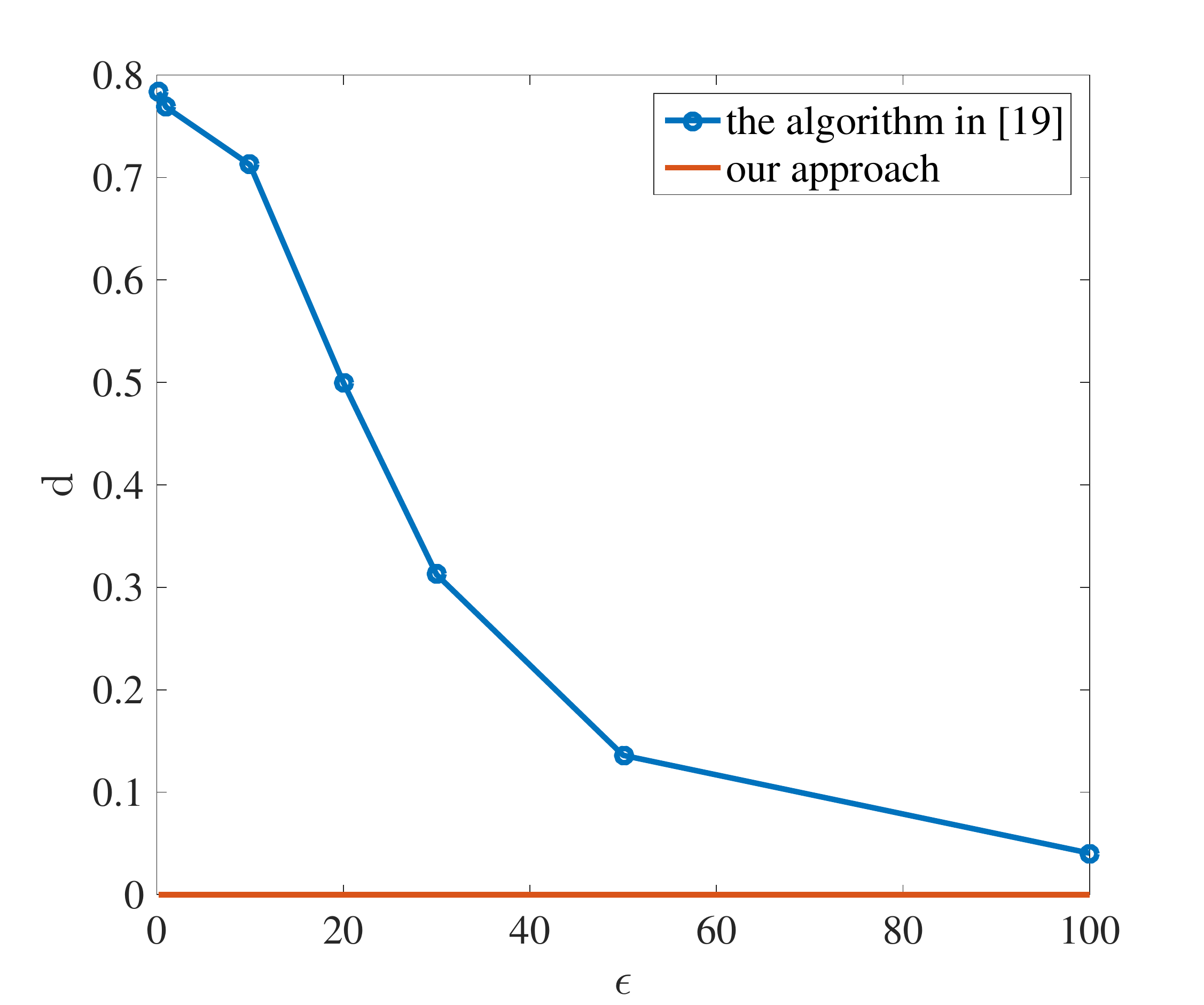}
    	\end{center}
    	\caption{The comparison of our approach with the algorithm in  \cite{huang2015differentially} in terms of optimization error.}
    	\label{fig:d}
    \end{figure}

    Simulation results from 5,000 runs showed that our approach  converged to $[0.35;0.45]$ with an error $d=3.14 \times10^{-14}$, which is negligible compared with the simulation results under the algorithm in  \cite{huang2015differentially} (cf. Fig. \ref{fig:d}, where each differential privacy level was implemented for 5,000 times). The results confirm the trade-off between privacy and accuracy for differential-privacy based approaches and demonstrate the advantages of our approach in terms of optimization accuracy.    

    \subsection{Comparison with the algorithm in \cite{lou2017privacy} }
    We  also compared our approach with the privacy-preserving optimization algorithm in \cite{lou2017privacy}. The network communication topology used for comparison is still the one in Fig. \ref{fig:network structure illustration} and the global objective function used is \eqref{eq:simulation} with $p_i$ $(i=1,2,..,6)$ fixed to $2$, $h_i$ $(i=1,2,..,6)$ fixed to $1$, and $\theta_i\in\mathbb{R}^2$. The adjacency matrix of network graph is defined in \eqref{eq: aij} for the algorithm in \cite{lou2017privacy}. Moreover, we let every agent update at each iteration and $c_i=1$ $(i=1,...,6)$ for  \cite{lou2017privacy}. The initial states are set to the same values for both algorithms.
    
\begin{figure}[!htbp]
    	\begin{center}
    		\includegraphics[width=0.49\textwidth]{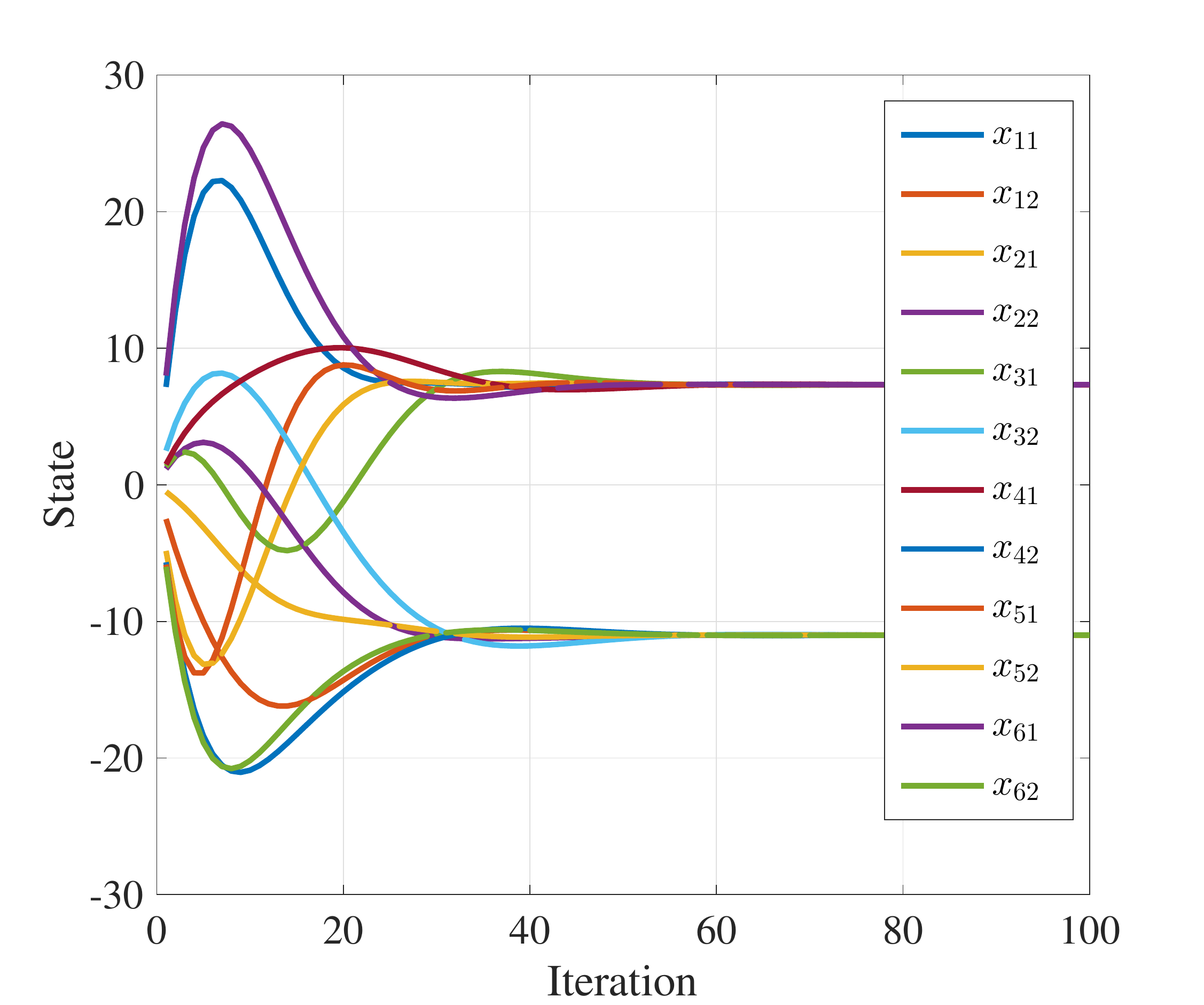}
    	\end{center}
    	\caption{The evolution of $\bm{x}_i$ in our approach.}
    	\label{fig:admm}
    \end{figure}

\begin{figure}[!htbp]
        	\begin{center}
        		\includegraphics[width=0.49\textwidth]{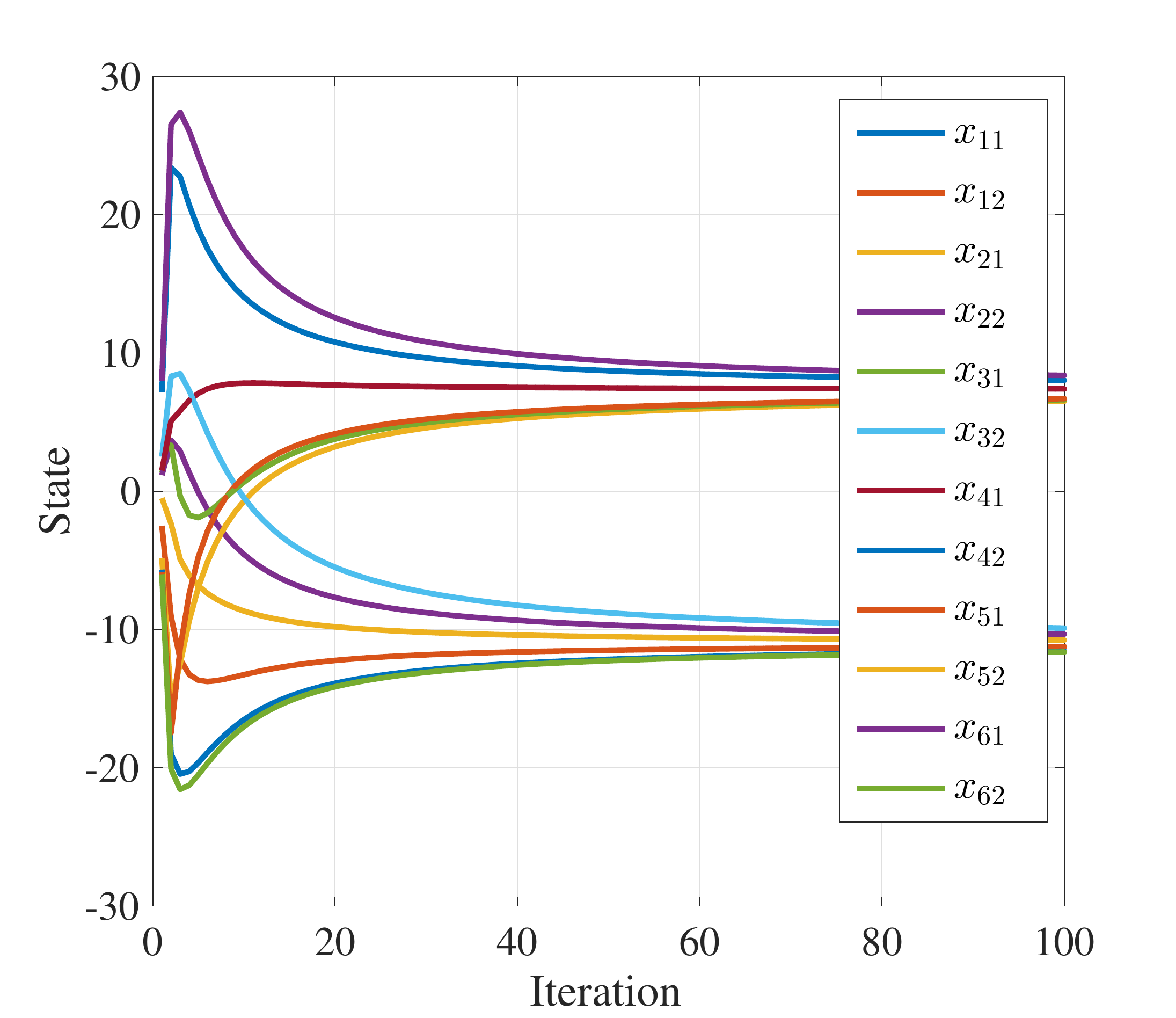}
        	\end{center}
        	\caption{The evolution of $\bm{x}_i$ in the algorithm of \cite{lou2017privacy}.}
        	\label{fig:subb}
        \end{figure}

     Fig. \ref{fig:admm} and Fig. \ref{fig:subb} show the evolution of $\bm{x}_i$ in our approach and the algorithm in \cite{lou2017privacy} respectively. It is clear that our approach converged faster than  the algorithm in \cite{lou2017privacy}.
    
        
\section{Implementation on Raspberry PI boards}
\begin{figure}[!htbp]
 	\begin{center}
 		\includegraphics[width=0.4\textwidth]{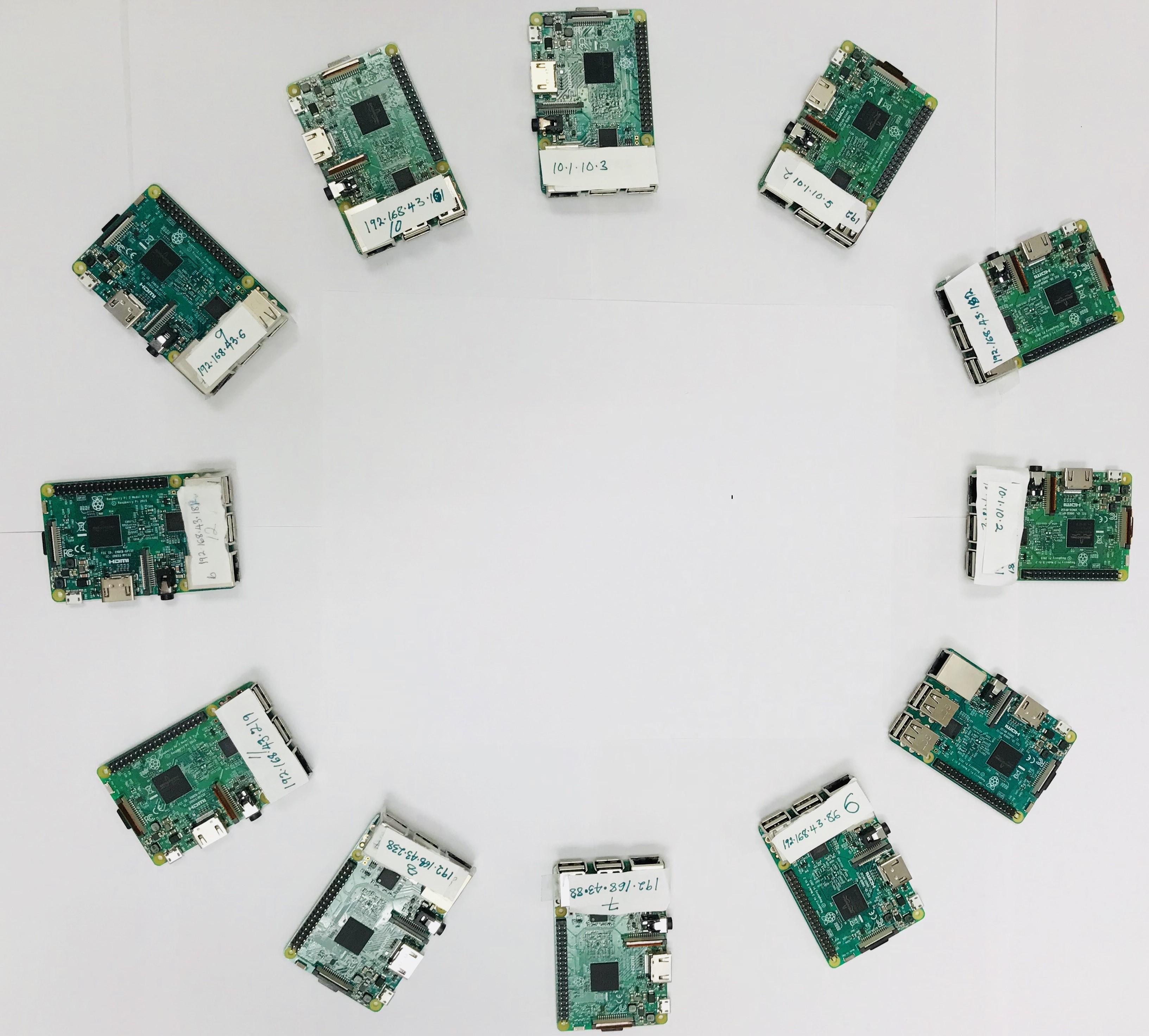}
 	\end{center}
 	\caption{The twelve Raspberry Pi boards}
 	\label{fig:experiment 1}
 \end{figure}

 We also implemented our privacy-preserving approach on twelve Raspberry Pi boards to confirm the efficiency of the
approach in real-world physical systems. Each board has 64-bit ARMv8
CPU and 1 GB RAM (cf. Fig. \ref{fig:experiment 1}). The optimization problem  \eqref{eq:simulation 1} was used in implementation with $p_i$ $(i=1,2,..,6)$ fixed to $2$, $h_i$ $(i=1,2,..,6)$ fixed to $1$, and $\theta_i\in\mathbb{R}$. In the implementation, ``libpaillier-0.8" library \cite{Paillier_Library} was used to realize the Paillier encryption and decryption process,  ``sys/socket.h" C library was used to conduct communication through Wi-Fi, and ``pthread" C library was
used to generate multiple parallel threads to realize parallelism in multi-agent networks.  
 The encryption and decryption keys were chosen as 512-bit long. 
 
\begin{figure}[!htbp]
        	\begin{center}
        		\includegraphics[width=0.49\textwidth]{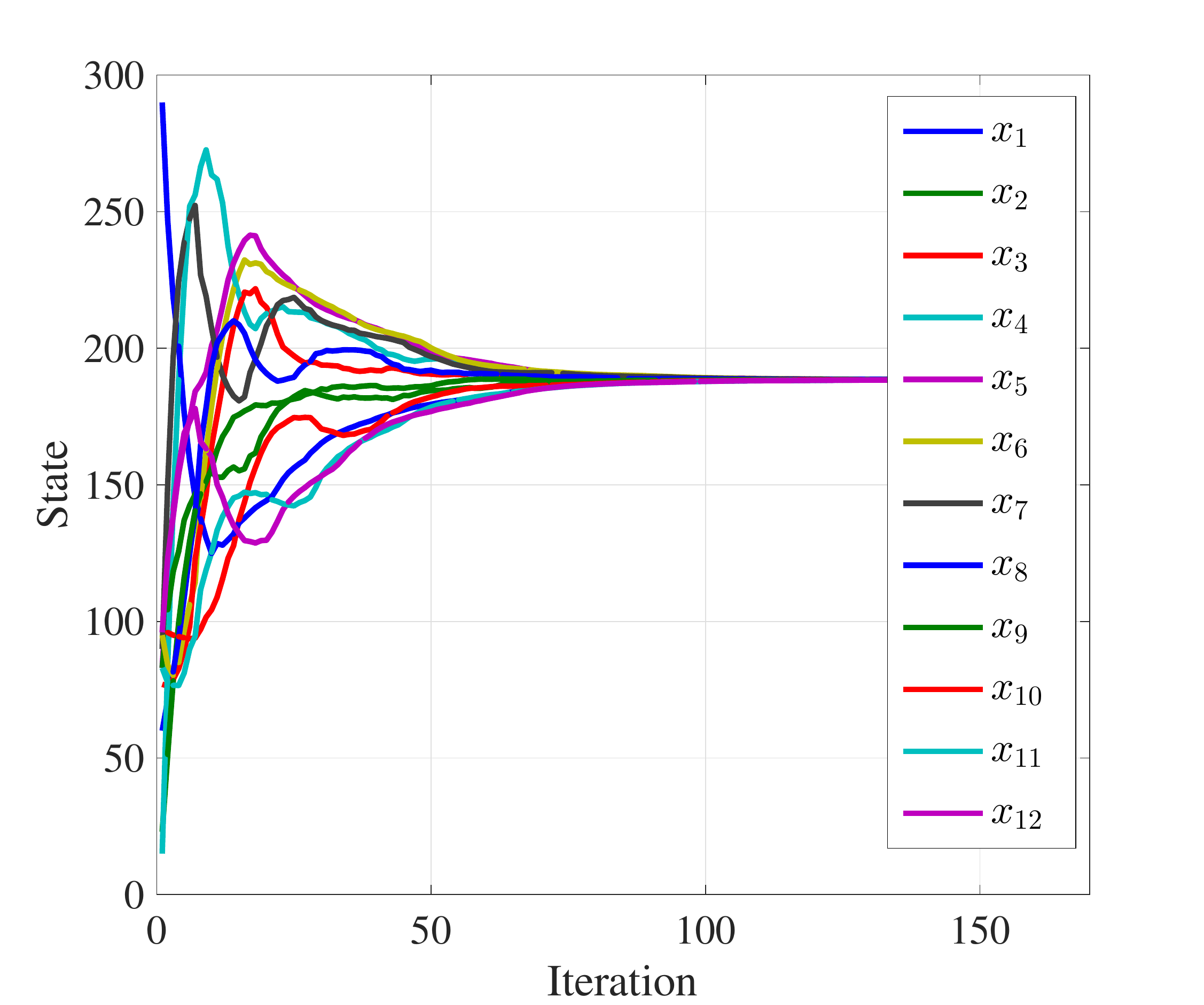}
        	\end{center}
        	\caption{The evolution of $x_i$ in the experimental verification using Raspberry Pi boards.}
        	\label{fig:experiment}
        \end{figure}

Implementation results confirmed that our approach always converged to the optimal solution. Fig. \ref{fig:experiment} visualizes the  evolution of $x_i$ $(i=1,2,...,12)$ in one specific  implementation where the network topology used is a cycle graph. We can see that each $x_i$ converged to the optimal solution  $188.417$.

\section{Conclusions}
In this paper, we presented a privacy-preserving decentralized optimization approach by proposing a new ADMM and leveraging partially homomorphic cryptography. By incorporating Paillier cryptosystem into the newly proposed decentralized ADMM, our approach provides guarantee for privacy preservation without compromising the solution in the absence of any aggregator or third party. This is in sharp contrast to differential-privacy based approaches which protect privacy through injecting noise and are subject to a fundamental trade-off between privacy and accuracy.  Theoretical analysis confirms that an honest-but-curious adversary cannot infer the information of neighboring agents even by recording and analyzing the information exchanged in multiple iterations.  The new ADMM allows time-varying penalty matrices and have a theoretically guaranteed convergence rate of $O(1/t)$, which makes it of mathematical interest by itself. Numerical and experimental results are given to confirm the effectiveness and efficiency of the proposed approach.




\section*{APPENDIX }
\subsection{Proof of Theorem \ref{Theorem:J-ADMM Convergence}}
The key idea to prove Theorem \ref{Theorem:J-ADMM Convergence} is to show that Algorithm  \uppercase\expandafter{\romannumeral1} converges to the saddle point of the Lagrangian function $L(\bm{x},\bm{\lambda})=f(\bm{x})+\bm{\lambda}^TA\bm{x}$. To achieve this goal, we introduce a variational inequality $MVI(Q,U)$ first and prove that the solution of  $MVI(Q,U)$  is also the saddle point of the Lagrangian function $L(\bm{x},\bm{\lambda})=f(\bm{x})+\bm{\lambda}^TA\bm{x}$ (which is formulated as Lemma 1). Then we introduce a sufficient condition for solving $MVI(Q,U)$ in Lemma 2. After the two steps, what  is left is to prove that the iterates of Algorithm \uppercase\expandafter{\romannumeral1} satisfy the condition in Lemma 2 when $k\to\infty$, i.e., Algorithm \uppercase\expandafter{\romannumeral1}  converges to the solution of  $MVI(Q,U)$ (Theorem \ref{Theorem:L} and Theorem \ref{Thereom:lim0}). 

We form a variational inequality $MVI(Q,U)$ similar to (5)-(6) in \cite{han2012note} first:
\begin{equation}
\begin{aligned}
\langle\bm{u}-\bm{u}^*,\bm{Q}(\bm{u}^*)\rangle\ge\bm{0}, \quad \forall \bm{u}, \label{eq:mvi 1}
\end{aligned}
\end{equation}
where
\begin{equation}  \label{eq:mvi 2}
\begin{aligned}
&{ \bm{u}^*:=\left(\begin{array}{c}
	\bm{x}_1^*\\
	\bm{x}_2^*\\
	\vdots\\
	\bm{x}_N^*\\
	\bm{\lambda}^*
	\end{array}\right),}
\quad
{ \bm{Q}(\bm{u}^*):=\left(\begin{array}{c}
	\xi_1^*+[A]_1^T\bm{\lambda}^*\\
	\xi_2^*+[A]_2^T\bm{\lambda}^*\\
	\vdots\\
	\xi_N^*+[A]_N^T\bm{\lambda}^*\\
	A\bm{x}^*
	\end{array}\right),}\\
& \xi_i^*\in\partial f_i(\bm{x}_i^*), \forall i\in\{1,2,...,N\}.
\end{aligned}
\end{equation}
In \eqref{eq:mvi 2}, $[A]_i$ denotes the columns of matrix $A$ that are associated with agent $i$. By recalling the first-order necessary and sufficient condition for convex programming \cite{han2012note}, it is easy to see that solving problem \eqref{eq:compact distributed admm form} amounts to solving the above $MVI(Q,U)$ \cite{han2012note}.  Denote the solution set of $MVI(Q,U)$ as $\mathcal{U}^*$. Since  $f_i$ is convex, $\partial f_i(\bm{x}_i)$ is monotone, the $MVI(Q,U)$ is solvable and $\mathcal{U}^*$ is nonempty \cite{han2012note}. 

Next, we introduce several lemmas and theorems that contribute to the proof of Theorem \ref{Theorem:J-ADMM Convergence}.
\begin{Lemma 1} \label{Lemma:2}
	Each $\bm{u}^*=(\bm{x}^*,\bm{\lambda}^*)$ in $\mathcal{U}^*$ is also the saddle point of the Lagrangian function $L(\bm{x},\bm{\lambda})=f(\bm{x})+\bm{\lambda}^TA\bm{x}$.
\end{Lemma 1}
{\it Proof}: The results can be obtained from Part 2.1 in \cite{he2016proximal} directly.  \hfill $\blacksquare$

\begin{Lemma 1} \label{Lemma lim}
	If $A\bm{x}^{k+1}=\bm{0}$ and $\bm{x}^{k+1}=\bm{x}^{k}$ hold, then $(\bm{x}_1^{k+1},\bm{x}_2^{k+1},...,\bm{x}_N^{k+1},\bm{\lambda}^{k+1})$ is a solution to $MVI(Q,U)$.
\end{Lemma 1}
{\it Proof}:  Using the definition of matrix A and the update rule of $\bm{\lambda}^{k+1}$ in \eqref{eq:time-varing admm 3}, we can see that the assumption  $A\bm{x}^{k+1}=\bm{0}$ implies $\bm{\lambda}^{k+1}=\bm{\lambda}^{k}$ and $\bm{x}_1^{k+1}=\bm{x}_2^{k+1}=...=\bm{x}_N^{k+1}$. 

On the other hand, we know that $\bm{x}_i^{k+1}$ is the optimizer of \eqref{eq:j-ADMM x update}.  By using the first-order optimality condition, we get
\begin{eqnarray}
\lefteqn{(\bm{x}_i-\bm{x}_i^{k+1})^T(\xi_i^{k+1}+\sum\limits_{j\in{\mathcal{N}_i}}(\bm{\lambda}_{i,j}^k}\nonumber \\
&  +\rho_{i,j}^k( \bm{x}_i^{k+1}-\bm{x}_{j}^k ))+ \gamma_i(\bm{x}_i^{k+1}-\bm{x}_i^k))\ge 0.
\end{eqnarray}
where $\xi_i^{k+1}\in  \partial f_i(\bm{x}_i^{k+1})$. Then based on the assumption $\bm{x}^{k+1}=\bm{x}^{k}$, the fact $\bm{\lambda}^{k+1}=\bm{\lambda}^{k}$, and the definition of matrix A, we have
$(\bm{x}_i-\bm{x}_i^{k+1})^T(\xi_i^{k+1}+[A]_i^T\bm{\lambda}^{k+1})\ge 0.$
Therefore, $(\bm{x}_1^{k+1},\bm{x}_2^{k+1},...,\bm{x}_N^{k+1},\bm{\lambda}^{k+1})$ is a solution to $MVI(Q,U)$.  \hfill $\blacksquare$

Lemma 2 provides a sufficient condition for solving $MVI(Q,U)$. According to Lemma 1, we know that the solution to $MVI(Q,U)$  is also the saddle point of the Lagrangian function. Next, we prove that the iterates in Algorithm \uppercase\expandafter{\romannumeral1} satisfy $\lim\limits_{k\to\infty}A\bm{x}^{k+1}=\bm{0}$ and $\lim\limits_{k\to\infty}\bm{x}^{k+1}-\bm{x}^{k}=\bm{0}$, i.e., Algorithm \uppercase\expandafter{\romannumeral1}  converges to the solution to  $MVI(Q,U)$. To achieve this goal, we first establish the relationship \eqref{eq:theorem L} about iterates $k$ and $k+1$ in Theorem \ref{Theorem:L}, whose proof is mainly based on convex properties. Then based on the relationship, we further prove   $\lim\limits_{k\to\infty}A\bm{x}^{k+1}=\bm{0}$ and $\lim\limits_{k\to\infty}\bm{x}^{k+1}-\bm{x}^{k}=\bm{0}$ in Theorem \ref{Thereom:lim0}.
 
\begin{Theorem 1} \label{Theorem:L}
		Let $\bm{\rho}^k$ satisfy Condition A, $\bar{Q}\triangleq Q_P+Q_C^k$ satisfy Condition B, and $(\bm{x}^*,\bm{\lambda}^*)$ be the saddle point of the Lagrangian function $L(\bm{x},\bm{\lambda})=f(\bm{x})+\bm{\lambda}^TA\bm{x}$, then we have
		\begin{equation} \label{eq:theorem L}
		\begin{aligned}
		&\quad\parallel \bm{\lambda}^{k+1}-\bm{\lambda}^*\parallel_{(\bm{\rho}^{k+1})^{-1}}^2+\parallel \bm{x}^{k+1}-\bm{x}^*\parallel_{\bar{Q}}^2\\
		&\le 	\parallel \bm{\lambda}^{k}-\bm{\lambda}^*\parallel_{(\bm{\rho}^k)^{-1}}^2+\parallel \bm{x}^{k}-\bm{x}^*\parallel_{\bar{Q}}^2 \\
	&	-(	\parallel A\bm{x}^{k+1}\parallel_{\bm{\rho}^k}^2+\parallel \bm{x}^{k+1}-\bm{x}^k\parallel_{-A^T\bm{\rho}^kA+\bar{Q}}^2) \\
	& +\parallel A\bm{x}^{k+1} \parallel_{\bm{\rho}^{k+1}}^2-\parallel A\bm{x}^{k} \parallel_{\bm{\rho}^k}^2.
		\end{aligned}
		\end{equation}
\end{Theorem 1}
 	
To prove Theorem \ref{Theorem:L}, we first introduce two lemmas:
\begin{Lemma 1} \label{Lemma_3}
	Let $\bm{x}^k=[\bm{x}_1^{kT},\bm{x}_2^{kT},...,\bm{x}_N^{kT} ]^T$ and $\bm{\lambda}^k=[\bm{\lambda}_{i,j}^k]_{ij,e_{i,j}\in E}$ be the intermediate results of iteration $k$ in Algorithm \uppercase\expandafter{\romannumeral1}, then the following inequality holds for all
	$k$:
	\begin{equation}
		\begin{aligned}
			\lefteqn{ f(\bm{x})-f(\bm{x}^{k+1})+(\bm{x}-\bm{x}^{k+1})^TA^T\bm{\lambda}^{k}+(\bm{x}-\bm{x}^{k+1})^T} \label{eq:lemma 3} \\
			& \bm{\cdot} A^T\bm{\rho}^kA\bm{x}^{k}+(\bm{x}-\bm{x}^{k+1})^T\bar{Q}(\bm{x}^{k+1}-\bm{x}^{k})\ge 0,
		\end{aligned}
	\end{equation}
	where  $\bar{Q}\triangleq Q_P+Q_C^k$.
\end{Lemma 1}

{\it Proof}: The proof follows from \cite{zhang2017distributed}. For completeness, we sketch the proof here. Denote by $g_i$ the function
\begin{equation}\label{eq:gi J-ADMM}
	\begin{aligned}
		g_i^k(\bm{x}_i)=\sum\limits_{j\in\mathcal{N}_{i}}(\bm{\lambda}_{i,j}^{kT}\bm{x}_{i}
		+ \frac{\rho_{i,j}^k}{2}\parallel \bm{x}_{i}-\bm{x}_{j}^{k} \parallel ^2)+\frac{\gamma_i}{2}\parallel \bm{x}_i-\bm{x}_i^k \parallel^2.
	\end{aligned}
\end{equation}

Using $\xi_i^{k+1}\in \partial f_i(\bm{x}_i^{k+1})$, we can get
 $\xi_i^{k+1}+\triangledown g_i(\bm{x}_i^{k+1})=\bm{0}$
and $(\bm{x}_i-\bm{x}_i^{k+1})^T[\xi_i^{k+1}+\triangledown g_i(\bm{x}_i^{k+1})]=0$ based on the fact that $\bm{x}_i^{k+1}$ is the optimizer
of $g_i^k+f_i$.
On the other hand, as $f_i$ is convex, the following relationship holds:
	$$f_i(\bm{x}_i)\ge f_i(\bm{x}_i^{k+1})+(\bm{x}_i-\bm{x}_i^{k+1})^T\xi_i^{k+1}.$$

Then we can get
	$f_i(\bm{x}_i)- f_i(\bm{x}_i^{k+1})+(\bm{x}_i-\bm{x}_i^{k+1})^T\triangledown g_i(\bm{x}_i^{k+1})\ge 0.$

Substituting $\triangledown g_i(\bm{x}_i^{k+1})$ with
\eqref{eq:gi J-ADMM}, we obtain
\begin{equation}
	\begin{aligned}
		\lefteqn{ f_i(\bm{x}_i)- f_i(\bm{x}_i^{k+1})+(\bm{x}_i-\bm{x}_i^{k+1})^T\bm{\cdot} }\nonumber \\
		& (\sum\limits_{j\in{\mathcal{N}_i}}(\bm{\lambda}_{i,j}^k+ \rho_{i,j}^k( \bm{x}_i^{k+1}-\bm{x}_{j}^k ))+ \gamma_i(\bm{x}_i^{k+1}-\bm{x}_i^k))\ge 0. \nonumber
	\end{aligned}
\end{equation}

Noting $\bm{\lambda}_{i,i}=\bm{0}$ and $\bm{\lambda}_{i,j}=-\bm{\lambda}_{j,i}$, based on the definition of matrices $A$ and $\bm{\rho}$, we can rewrite the above inequality as
\begin{equation}
	\begin{aligned}
		\lefteqn{ f_i(\bm{x}_i)- f_i(\bm{x}_i^{k+1})+(\bm{x}_i-\bm{x}_i^{k+1})^T\bm{\cdot}} \label{eq:2} \\
		&  ([A]_i^T\bm{\lambda}^{k} +\sum\limits_{j\in{\mathcal{N}_i}} \rho_{i,j}^k ( \bm{x}_i^{k+1}-\bm{x}_{j}^k )+\gamma_i(\bm{x}_i^{k+1}-\bm{x}_i^k))\ge 0.
	\end{aligned}
\end{equation}

Summing both sides of \eqref{eq:2} over $i = 1,2, \ldots,N$, and using
\begin{equation}
	\begin{aligned}	
		\lefteqn{\sum\limits_{i=1}^{N}(\bm{x}_i-\bm{x}_i^{k+1})^T[ A]_i^T\bm{\lambda}^{k}=(\bm{x}-\bm{x}^{k+1})^TA^T\bm{\lambda}^{k},} \nonumber\\	
			&\sum\limits_{i=1}^{N}(\bm{x}_i-\bm{x}_i^{k+1})^T 
			\sum\limits_{j\in{\mathcal{N}_i}}\rho_{i,j}^k\bm{x}_i^{k+1} 
			=(\bm{x}-\bm{x}^{k+1})^TQ_C^k\bm{x}^{k+1}, \nonumber \\	
			&\sum\limits_{i=1}^{N}(\bm{x}_i-\bm{x}_i^{k+1})^T
			 \sum\limits_{j\in{\mathcal{N}_i}} \rho_{i,j}^k\bm{x}_{j}^k \nonumber \\
		&\qquad\qquad\qquad=(\bm{x}-\bm{x}^{k+1})^T(-A^T\bm{\rho}^kA+Q_C^k)\bm{x}^{k}, \nonumber\\			
			&\sum\limits_{i=1}^{N}(\bm{x}_i-\bm{x}_i^{k+1})^T 
			\gamma_i(\bm{x}_i^{k+1}-\bm{x}_i^{k}) \nonumber \\
			& \qquad\qquad\qquad=(\bm{x}-\bm{x}^{k+1})^TQ_P(\bm{x}^{k+1}-\bm{x}^{k}), \nonumber 
	\end{aligned}
\end{equation}

we can get the lemma. \hfill $\blacksquare$

\begin{Lemma 2} \label{Lemma_4}
		 Let $\bm{x}^k=[\bm{x}_1^{kT},\bm{x}_2^{kT},...,\bm{x}_N^{kT} ]^T$ and
	$\bm{\lambda}^k=[\bm{\lambda}_{i,j}^k]_{ij,e_{i,j}\in E}$ be
	the intermediate results of iteration $k$ in
	 Algorithm \uppercase\expandafter{\romannumeral1}, then the following equality holds for all
	$k$:
	\begin{equation}
		\begin{aligned}
			&-(\bm{x}^{k+1})^TA^T(\bm{\lambda}^{k}-\bm{\lambda}^{*}) \label{eq:lemma 4} \\
			& -(\bm{x}^{k+1})^TA^T\bm{\rho}^kA\bm{x}^{k}+(\bm{x}^{*}-\bm{x}^{k+1})^T\bar{Q}(\bm{x}^{k+1}-\bm{x}^{k})  \\
			&=-\frac{1}{2}(\parallel \bm{\lambda}^{k+1}-\bm{\lambda}^*\parallel_{(\bm{\rho}^{k+1})^{-1}}^2-\parallel \bm{\lambda}^{k}-\bm{\lambda}^*\parallel_{(\bm{\rho}^{k+1})^{-1}}^2)  \\
			&+\frac{1}{2}\parallel \bm{\lambda}^{k+1}-\bm{\lambda}^k\parallel_{(\bm{\rho}^{k+1})^{-1}}^2+\frac{1}{2}\parallel A(\bm{x}^{k+1}-\bm{x}^k)\parallel_{\bm{\rho}^k}^2  \\
			& -\frac{1}{2}\parallel A\bm{x}^{k+1}\parallel_{\bm{\rho}^k}^2-\frac{1}{2}\parallel A\bm{x}^{k}\parallel_{\bm{\rho}^k}^2 -\frac{1}{2}\parallel \bm{x}^{k+1}-\bm{x}^k \parallel_{\bar{Q}}^2\\
			& -\frac{1}{2}(\parallel \bm{x}^{k+1}-\bm{x}^*\parallel_{\bar{Q}}^2-\parallel \bm{x}^{k}-\bm{x}^*\parallel_{\bar{Q}}^2). 
		\end{aligned}
	\end{equation}
\end{Lemma 2}

{\it Proof}:  For a scalar $a$, we have $a^T=a$. Recall
$\bm{\lambda}^{k+1}=\bm{\lambda}^k+\bm{\rho}^{k+1}A\bm{x}^{k+1}$ and notice that $\bm{\rho}^{k+1}$ is a positive definite diagonal matrix, we can get
\begin{equation}
	\begin{aligned}
		(\bm{x}^{k+1})^TA^T(\bm{\lambda}^{k}-\bm{\lambda}^{*})
		=(\bm{\lambda}^{k+1}-\bm{\lambda}^k)^T(\bm{\rho}^{k+1})^{-1}(\bm{\lambda}^{k}-\bm{\lambda}^{*}). \label{eq:lemma2 1}
	\end{aligned}
\end{equation}

On the other hand, since $(\bm{x}^*,\bm{\lambda}^*)$ is the saddle point of the Lagrangian function \eqref{eq:lagrangian function}, we can get
$A\bm{x}^*=\bm{0}$ \cite{wei2012distributed}. Moreover, the following
equalities can be established by using algebraic manipulations:
\begin{equation}
\begin{aligned}
\lefteqn{(\bm{x}^{k+1}-\bm{x}^{*})^T\bar{Q}(\bm{x}^{k+1}-\bm{x}^{k})=\frac{1}{2}\parallel \bm{x}^{k+1}-\bm{x}^{k}\parallel_{\bar{Q}}^2} \\
&\qquad\qquad\quad+\frac{1}{2}(\parallel \bm{x}^{k+1}-\bm{x}^*\parallel_{\bar{Q}}^2-\parallel \bm{x}^{k}-\bm{x}^*\parallel_{\bar{Q}}^2),  \label{eq:lemma2 2} 
\end{aligned}
\end{equation}
\begin{equation}
\begin{aligned}
\lefteqn{-\bm{x}^{(k+1)T}A^T\bm{\rho}^kA\bm{x}^{k}=\frac{1}{2}\parallel A(\bm{x}^{k+1}-\bm{x}^{k})\parallel_{\bm{\rho}^k}^2} \\
&\qquad\qquad\qquad-\frac{1}{2}\parallel A \bm{x}^{k+1}\parallel_{\bm{\rho}^k}^2-\frac{1}{2}\parallel A\bm{x}^{k}\parallel_{\bm{\rho}^k}^2,  \label{eq:lemma2 3} 
\end{aligned}
\end{equation}
\begin{equation}
	\begin{aligned}
		\lefteqn{(\bm{\lambda}^{k+1}-\bm{\lambda}^k)^T(\bm{\rho}^{k+1})^{-1}(\bm{\lambda}^{k}-\bm{\lambda}^{*})} \\
		&=\frac{1}{2}(\parallel \bm{\lambda}^{k+1}-\bm{\lambda}^{*}\parallel_{(\bm{\rho}^{k+1})^{-1}}^2-\parallel \bm{\lambda}^{k}-\bm{\lambda}^{*}\parallel_{(\bm{\rho}^{k+1})^{-1}}^2), \\
		&-\frac{1}{2}\parallel \bm{\lambda}^{k+1}-\bm{\lambda}^{k}\parallel_{(\bm{\rho}^{k+1})^{-1}}^2.
\label{eq:lemma2 4}
	\end{aligned}
\end{equation}
Then we can obtain \eqref{eq:lemma 4} by plugging equalities \eqref{eq:lemma2 1}-\eqref{eq:lemma2 4} 
into the left hand side of \eqref{eq:lemma 4}. \hfill $\blacksquare$

Now we can proceed to prove Theorem \ref{Theorem:L}. By setting $\bm{x}=\bm{x}^*$ in \eqref{eq:lemma 3}, we can get
\begin{equation}
	\begin{aligned}
		\lefteqn{ f(\bm{x}^*)-f(\bm{x}^{k+1})+(\bm{x}^*-\bm{x}^{k+1})^TA^T\bm{\lambda}^{k}+(\bm{x}^*-\bm{x}^{k+1})^T} \nonumber \\
		&\qquad\bm{\cdot}A^T\bm{\rho}^kA\bm{x}^{k}+(\bm{x}^*-\bm{x}^{k+1})^T\bar{Q}(\bm{x}^{k+1}-\bm{x}^{k})\ge 0. \nonumber
	\end{aligned}
\end{equation}
Recalling $A\bm{x}^*=\bm{0}$, the above inequality can be rewritten as
\begin{equation}
	\begin{aligned}
		\lefteqn{ f(\bm{x}^*)-f(\bm{x}^{k+1})-\bm{x}^{(k+1)T}A^T\bm{\lambda}^{k} } \label{eq:1 for theorem 2} \\
		&-\bm{x}^{(k+1)T}A^T\bm{\rho}^kA\bm{x}^{k} +(\bm{x}^*-\bm{x}^{k+1})^T\bar{Q}(\bm{x}^{k+1}-\bm{x}^{k})\ge 0. 
	\end{aligned}
\end{equation}

Now adding and subtracting the term $\bm{\lambda}^{*T}A\bm{x}^{k+1}$ from the left hand side of \eqref{eq:1 for theorem 2} gives
\begin{equation} \label{eq:theorem 2}
	\begin{aligned}
	\lefteqn{ f(\bm{x}^*)-f(\bm{x}^{k+1})-\bm{\lambda}^{*T}A\bm{x}^{k+1}-\bm{x}^{(k+1)T}A^T(\bm{\lambda}^{k}-\bm{\lambda}^{*}) }\\
	&-\bm{x}^{(k+1)T}A^T\bm{\rho}^kA\bm{x}^{k} +(\bm{x}^*-\bm{x}^{k+1})^T\bar{Q}(\bm{x}^{k+1}-\bm{x}^{k})\ge 0. 
	\end{aligned}
\end{equation}

Using $L(\bm{x},\bm{\lambda}^*)-L(\bm{x}^*,\bm{\lambda}^*)\ge 0$ and $A\bm{x}^*=\bm{0}$,  we have
\begin{equation} \nonumber
\begin{aligned}
 \lefteqn{-\bm{x}^{(k+1)T}A^T(\bm{\lambda}^{k}-\bm{\lambda}^{*})}\\
 &-\bm{x}^{(k+1)T}A^T\bm{\rho}^kA\bm{x}^{k} +(\bm{x}^*-\bm{x}^{k+1})^T\bar{Q}(\bm{x}^{k+1}-\bm{x}^{k}) \\
 &\ge f(\bm{x}^{k+1})+\bm{\lambda}^{*T}A\bm{x}^{k+1}-f(\bm{x}^*) \ge 0. 
\end{aligned}
\end{equation}

Now by plugging \eqref{eq:lemma 4} into the left hand side of the above inequality, we can obtain
\begin{eqnarray}
	&& -\frac{1}{2}(\parallel \bm{\lambda}^{k+1}-\bm{\lambda}^{*}\parallel_{(\bm{\rho}^{k+1})^{-1}}^2-\parallel \bm{\lambda}^{k}-\bm{\lambda}^{*}\parallel_{(\bm{\rho}^{k+1})^{-1}}^2) \nonumber \\
	&&+\frac{1}{2}\parallel \bm{\lambda}^{k+1}-\bm{\lambda}^k\parallel_{(\bm{\rho}^{k+1})^{-1}}^2+\frac{1}{2}\parallel A\bm{x}^{k+1}-A\bm{x}^{k} \parallel_{\bm{\rho}^k}^2 \nonumber \\
	&&-\frac{1}{2}\parallel A\bm{x}^{k+1} \parallel_{\bm{\rho}^k}^2-\frac{1}{2}\parallel A\bm{x}^{k} \parallel_{\bm{\rho}^k}^2
  -\frac{1}{2}\parallel \bm{x}^{k+1}-\bm{x}^*\parallel_{\bar{Q}}^2\nonumber \\
	&&+\frac{1}{2}\parallel \bm{x}^{k}-\bm{x}^*\parallel_{\bar{Q}}^2-
	\frac{1}{2}\parallel \bm{x}^{k+1}-\bm{x}^k \parallel_{\bar{Q}}^2\ge 0. \nonumber
\end{eqnarray}

Noting $\parallel \bm{\lambda}^{k+1}-\bm{\lambda}^k\parallel_{(\bm{\rho}^{k+1})^{-1}}^2=\parallel A\bm{x}^{k+1} \parallel_{\bm{\rho}^{k+1}}^2$, the above inequality can be rewritten as

	\begin{equation}
	\begin{aligned}
	&\quad\parallel \bm{\lambda}^{k+1}-\bm{\lambda}^*\parallel_{(\bm{\rho}^{k+1})^{-1}}^2+\parallel \bm{x}^{k+1}-\bm{x}^*\parallel_{\bar{Q}}^2\\
	&\le 	\parallel \bm{\lambda}^{k}-\bm{\lambda}^*\parallel_{(\bm{\rho}^{k+1})^{-1}}^2+\parallel \bm{x}^{k}-\bm{x}^*\parallel_{\bar{Q}}^2 \\
	&	-(	\parallel A\bm{x}^{k+1}\parallel_{\bm{\rho}^k}^2+\parallel \bm{x}^{k+1}-\bm{x}^k\parallel_{-A^T\bm{\rho}^kA+\bar{Q}}^2) \\
	&+\parallel A\bm{x}^{k+1} \parallel_{\bm{\rho}^{k+1}}^2-\parallel A\bm{x}^{k} \parallel_{\bm{\rho}^k}^2.
	\end{aligned}
	\end{equation}

Recall that from Condition A, $\bm{\rho}^{k+1}\succeq\bm{\rho}^k$ and $\bm{\rho}^k$ $(k=1,2,...)$ are positive definite diagonal matrices. So we have  $(\bm{\rho}^{k+1})^{-1}\preceq(\bm{\rho}^k)^{-1}$ \cite{kontogiorgis1998variable}, and consequently $\parallel \bm{\lambda}^{k}-\bm{\lambda}^*\parallel_{(\bm{\rho}^{k+1})^{-1}}^2\preceq\parallel \bm{\lambda}^{k}-\bm{\lambda}^*\parallel_{(\bm{\rho}^{k})^{-1}}^2$, which proves Theorem \ref{Theorem:L}. 

\hfill $\blacksquare$

Theorem \ref{Theorem:L} established the relationship between iterates $k$ and $k+1$ in Algorithm \uppercase\expandafter{\romannumeral1}. Based on this relationship, we can have the following theorem which shows that Algorithm \uppercase\expandafter{\romannumeral1}  converges to the solution to  $MVI(Q,U)$.

 \begin{Theorem 1} \label{Thereom:lim0}
 	Let $\bm{u}^{k}=(\bm{x}^{k},\bm{\lambda}^{k})$ be the sequence generated by Algorithm \uppercase\expandafter{\romannumeral1}, then we have
 	\begin{eqnarray}
 	\lim\limits_{k\to\infty}(\parallel A\bm{x}^{k+1}\parallel_{\bm{\rho}^k}^2+\parallel \bm{x}^{k+1}-\bm{x}^k\parallel_{-A^T\bm{\rho}^kA+\bar{Q}}^2)=0.
 	\end{eqnarray}
 \end{Theorem 1}
 {\it Proof:} Let 
 $\alpha^k=\parallel \bm{\lambda}^{k}-\bm{\lambda}^*\parallel_{(\bm{\rho}^k)^{-1}}^2+\parallel \bm{x}^{k}-\bm{x}^*\parallel_{\bar{Q}}^2$. According to Theorem \ref{Theorem:L}, we have
 	\begin{equation}
 	\begin{aligned}
 	&\alpha^{k+1}\le \alpha^{k}	 +\parallel A\bm{x}^{k+1} \parallel_{\bm{\rho}^{k+1}}^2-\parallel A\bm{x}^{k} \parallel_{\bm{\rho}^k}^2 \\
 	&	-(	\parallel A\bm{x}^{k+1}\parallel_{\bm{\rho}^k}^2+\parallel \bm{x}^{k+1}-\bm{x}^k\parallel_{-A^T\bm{\rho}^kA+\bar{Q}}^2) \\
 	& \le ...\\
 	& \le   \alpha^{0}+\parallel A\bm{x}^{k+1} \parallel_{\bm{\rho}^{k+1}}^2-\parallel A\bm{x}^{0} \parallel_{\bm{\rho}^0}^2 \\
 	& 	-\sum\limits_{i=0}^{k}(	\parallel A\bm{x}^{i+1}\parallel_{\bm{\rho}^i}^2+\parallel \bm{x}^{i+1}-\bm{x}^i\parallel_{-A^T\bm{\rho}^iA+\bar{Q}}^2)\\
 	& \le   \alpha^{0}+\parallel \bm{x}^{k+1}-\bm{x}^* \parallel_{A^T\bm{\rho}^{k+1}A}^2\\
 	& 	-\sum\limits_{i=0}^{k}(	\parallel A\bm{x}^{i+1}\parallel_{\bm{\rho}^i}^2+\parallel \bm{x}^{i+1}-\bm{x}^i\parallel_{-A^T\bm{\rho}^iA+\bar{Q}}^2).
 	\end{aligned}
 	\end{equation} 	
 The last inequality comes from the fact that $A\bm{x}^*=\bm{0}$ and $\parallel A\bm{x}^{k+1}-A\bm{x}^* \parallel_{\bm{\rho}^{k+1}}^2$ can be written as $\parallel \bm{x}^{k+1}-\bm{x}^* \parallel_{A^T\bm{\rho}^{k+1}A}^2$. Recall that $\bm{\rho}^{0}\preceq\bm{\rho}^{k}\preceq\bm{\rho}^{k+1}\preceq\bar{\bm{\rho}}$ holds and $\bar{Q}-A^T\bar{\bm{\rho}}A$ is positive definite. Moving the term $\parallel \bm{x}^{k+1}-\bm{x}^* \parallel_{A^T\bm{\rho}^{k+1}A}^2$ to the left hand side of the above inequality, we have
 	 	\begin{equation}
 	 	\begin{aligned}
 &\quad	\lim\limits_{k\to\infty}( \alpha^{k+1}-\parallel \bm{x}^{k+1}-\bm{x}^* \parallel_{A^T\bm{\rho}^{k+1}A}^2 )\\
 &	=\lim\limits_{k\to\infty}(\parallel \bm{\lambda}^{k+1}-\bm{\lambda}^*\parallel_{(\bm{\rho}^{k+1})^{-1}}^2+ \parallel \bm{x}^{k+1}-\bm{x}^* \parallel_{\bar{Q}-A^T\bm{\rho}^{k+1}A}^2)\\
 &\ge 0
 	 	\end{aligned}
 	 	\end{equation} 
 Since $\alpha^0$ is positive and bounded and $\parallel A\bm{x}^{i+1}\parallel_{\bm{\rho}^i}^2+\parallel \bm{x}^{i+1}-\bm{x}^i\parallel_{-A^T\bm{\rho}^iA+\bar{Q}}^2$ is nonnegative,  following Theorem 3 in  \cite{he2002new}, we have
 \begin{eqnarray} \label{eq:limit 0}
 \lim\limits_{k\to\infty}(\parallel A\bm{x}^{k+1}\parallel_{\bm{\rho}^k}^2+\parallel \bm{x}^{k+1}-\bm{x}^k\parallel_{-A^T\bm{\rho}^kA+\bar{Q}}^2)=0. 
 \end{eqnarray} \hfill $\blacksquare$
  
Given that $\bm{\rho}^{k}$ satisfies Condition A and $\bar{Q}$ satisfies Condition B, we have that both $-A^T\bm{\rho}^kA+\bar{Q}$ and $\bm{\rho}^k$ are positive symmetric definite. Then according to Theorem \ref{Thereom:lim0}, we have $ A\bm{x}^{k+1}=\bm{0}$ and $\bm{x}^{k+1}=\bm{x}^k$ when $k\to \infty$.
 
 Therefore, based on Lemma \ref{Lemma lim}, we have that $(\bm{x}^{k+1},\bm{\lambda}^{k+1})$ in Algorithm \uppercase\expandafter{\romannumeral1} converges to a solution to $MVI(Q,U)$, i.e., a saddle point of the Lagrangian function \eqref{eq:lagrangian function} according to Lemma \ref{Lemma:2}. Since the objective function is convex, we can conclude Theorem \ref{Theorem:J-ADMM Convergence} \cite{wei2012distributed}. \hfill $\blacksquare$

\subsection{Proof of Theorem \ref{Theorem: convergence rate}}
 Now we prove that the convergence rate of Algorithm \uppercase\expandafter{\romannumeral1} is $O(1/t)$.
By plugging \eqref{eq:lemma 4} into the left hand side of \eqref{eq:theorem 2}, we can obtain
\begin{eqnarray}
&&f(\bm{x}^*)-f(\bm{x}^{k+1})-\bm{\lambda}^{*T}A\bm{x}^{k+1}\\
&& -\frac{1}{2}(\parallel \bm{\lambda}^{k+1}-\bm{\lambda}^{*}\parallel_{(\bm{\rho}^{k+1})^{-1}}^2-\parallel \bm{\lambda}^{k}-\bm{\lambda}^{*}\parallel_{(\bm{\rho}^{k+1})^{-1}}^2) \nonumber \\
&&+\frac{1}{2}\parallel \bm{\lambda}^{k+1}-\bm{\lambda}^k\parallel_{(\bm{\rho}^{k+1})^{-1}}^2+\frac{1}{2}\parallel A\bm{x}^{k+1}-A\bm{x}^{k} \parallel_{\bm{\rho}^k}^2 \nonumber \\
&&-\frac{1}{2}\parallel A\bm{x}^{k+1} \parallel_{\bm{\rho}^k}^2-\frac{1}{2}\parallel A\bm{x}^{k} \parallel_{\bm{\rho}^k}^2
-\frac{1}{2}\parallel \bm{x}^{k+1}-\bm{x}^*\parallel_{\bar{Q}}^2\nonumber \\
&&+\frac{1}{2}\parallel \bm{x}^{k}-\bm{x}^*\parallel_{\bar{Q}}^2-
\frac{1}{2}\parallel \bm{x}^{k+1}-\bm{x}^k \parallel_{\bar{Q}}^2\ge 0. \nonumber
\end{eqnarray}

Summing both sides of the above inequality over $k=0,1,...,t$, we have
\begin{eqnarray}
&&(t+1)f(\bm{x}^*)-\sum_{k=0}^{t}f(\bm{x}^{k+1})-\bm{\lambda}^{*T}A\sum_{k=0}^{t}\bm{x}^{k+1}  \nonumber\\
&& -\frac{1}{2}\parallel \bm{\lambda}^{t+1}-\bm{\lambda}^{*}\parallel_{(\bm{\rho}^{t+1})^{-1}}^2+\frac{1}{2}\parallel \bm{\lambda}^{0}-\bm{\lambda}^{*}\parallel_{(\bm{\rho}^{1})^{-1}}^2 \nonumber\\
&&-\sum_{k=1}^{t}\frac{1}{2}(\parallel \bm{\lambda}^{k}-\bm{\lambda}^{*}\parallel_{(\bm{\rho}^{k})^{-1}}^2-\parallel \bm{\lambda}^{k}-\bm{\lambda}^{*}\parallel_{(\bm{\rho}^{k+1})^{-1}}^2) \nonumber \\
&&+\frac{1}{2}\parallel A\bm{x}^{t+1}\parallel_{\bm{\rho}^{t+1}}^2
-\sum_{k=0}^{t}\frac{1}{2}\parallel A\bm{x}^{k+1} \parallel_{\bm{\rho}^k}^2-\frac{1}{2}\parallel A\bm{x}^{0} \parallel_{\bm{\rho}^0}^2\nonumber\\
&&-\frac{1}{2}\parallel \bm{x}^{t+1}-\bm{x}^*\parallel_{\bar{Q}}^2+\frac{1}{2}\parallel \bm{x}^{0}-\bm{x}^*\parallel_{\bar{Q}}^2\nonumber \\
&&-
\sum_{k=0}^{t}\frac{1}{2}\parallel \bm{x}^{k+1}-\bm{x}^k \parallel_{\bar{Q}-A^T\bm{\rho}^kA}^2\ge 0. \nonumber
\end{eqnarray}

 Following the above inequality, It is easy to obtain 
\begin{eqnarray}
&&(t+1)f(\bm{x}^*)-\sum_{k=0}^{t}f(\bm{x}^{k+1})-\bm{\lambda}^{*T}A\sum_{k=0}^{t}\bm{x}^{k+1}  \nonumber\\
&&+\frac{1}{2}\parallel \bm{\lambda}^{0}-\bm{\lambda}^{*}\parallel_{(\bm{\rho}^{1})^{-1}}^2 +\frac{1}{2}\parallel \bm{x}^{0}-\bm{x}^*\parallel_{\bar{Q}}^2\nonumber\\
&&+\frac{1}{2}\parallel A\bm{x}^{t+1}\parallel_{\bm{\rho}^{t+1}}^2
-\frac{1}{2}\parallel A\bm{x}^{t+1} \parallel_{\bm{\rho}^t}^2\ge 0. \nonumber
\end{eqnarray}

Recall that in \eqref{eq:limit 0}, we have proven  $\lim\limits_{k\to\infty}\parallel A\bm{x}^{k+1}\parallel_{\bm{\rho}^{k}}^2=0$. Then the relationship $\bm{\rho}^{0}\preceq\bm{\rho}^{k}\preceq\bm{\rho}^{k+1}\preceq\bar{\bm{\rho}}$ implies
 	$$\lim\limits_{t\to\infty}(\frac{1}{2}\parallel A\bm{x}^{t+1}\parallel_{\bm{\rho}^{t+1}}^2
 	-\frac{1}{2}\parallel A\bm{x}^{t+1} \parallel_{\bm{\rho}^t}^2)=0.$$
 Therefore,  there exists some constant $c$ such that $$\frac{1}{2}\parallel A\bm{x}^{t+1}\parallel_{\bm{\rho}^{t+1}}^2
-\frac{1}{2}\parallel A\bm{x}^{t+1} \parallel_{\bm{\rho}^t}^2\le c.$$

On the other hand, as our function is convex, we have $\sum_{k=0}^{t}f(\bm{x}^{k+1})\ge (t+1)f(\bar{\bm{x}}^{t+1})$ where $\bar{\bm{x}}^{t+1}=\frac{1}{t+1}\sum_{k=0}^{t}\bm{x}^{k+1}$. Therefore, we have
\begin{eqnarray}
&&(t+1)f(\bm{x}^*)-(t+1)f(\bar{\bm{x}}^{t+1})-(t+1)\bm{\lambda}^{*T}A\bar{\bm{x}}^{t+1}  \nonumber\\
&&+\frac{1}{2}\parallel \bm{\lambda}^{0}-\bm{\lambda}^{*}\parallel_{(\bm{\rho}^{1})^{-1}}^2 +\frac{1}{2}\parallel \bm{x}^{0}-\bm{x}^*\parallel_{\bar{Q}}^2+c\ge 0. \nonumber
\end{eqnarray}

 By dividing both sides by $-(t+1)$, we can obtain
\begin{eqnarray}
&&f(\bar{\bm{x}}^{t+1})+\bm{\lambda}^{*T}A\bar{\bm{x}}^{t+1}-f(\bm{x}^*)  \nonumber\\
&&\le \frac{1}{t+1}(\frac{1}{2}\parallel \bm{\lambda}^{0}-\bm{\lambda}^{*}\parallel_{(\bm{\rho}^{1})^{-1}}^2 +\frac{1}{2}\parallel \bm{x}^{0}-\bm{x}^*\parallel_{\bar{Q}}^2+c). \nonumber
\end{eqnarray}

Combining the above relationship with the Lagrangian function \eqref{eq:lagrangian function}, we can conclude Theorem \ref{Theorem: convergence rate}.



\bibliographystyle{unsrt}
\bibliography{abbr_bibli}


\end{document}